\documentclass[11pt,reqno]{amsart}      
\usepackage[english]{babel}
\usepackage[utf8x]{inputenc}
\usepackage[T1]{fontenc}
\usepackage[letterpaper, margin=1.2in]{geometry}
\usepackage{mathtools}
\usepackage[dvipsnames]{xcolor}
\usepackage{hyperref}
\usepackage{dsfont}
\usepackage{enumerate}
\usepackage{dsfont}
\usepackage[pdftex]{graphicx}
\usepackage{caption}
\usepackage{subcaption}
\usepackage{tikz}
\usetikzlibrary{positioning,decorations.pathreplacing,shapes.misc}
\tikzset{cross/.style={cross out, draw=black, fill=none, minimum size=2*(#1-\pgflinewidth), inner sep=0pt, outer sep=0pt}, cross/.default={2pt}}

\allowdisplaybreaks

\usepackage{amsmath, amsthm, amssymb}
\usepackage{latexsym,graphicx,multicol,mathrsfs,xypic}

\newtheorem{theorem}{Theorem}[section]
\newtheorem{definition}[theorem]{Definition}
\newtheorem{proposition}[theorem]{Proposition}
\newtheorem{lemma}[theorem]{Lemma}
\newtheorem{corollary}[theorem]{Corollary}

\newtheorem{remark}[theorem]{Remark}

\newtheorem*{proposition*}{Proposition}

\DeclareMathOperator{\var}{var}
\DeclareMathOperator{\cov}{cov}

\DeclareMathOperator{\supp}{supp}

\usepackage{autonum}

\begin{document}

\title[Equivalence of the grand canonical ensemble and the canonical ensemble]{Equivalence of the grand canonical ensemble and the canonical ensemble on 1d-lattice systems}

\author{Younghak Kwon}
\address{Department of Mathematics, University of California, Los Angeles}
\email{yhkwon@math.ucla.edu}

\author{Jaehun Lee}
\address{Department of Mathematics, Seoul National University}
\email{hun618@snu.ac.kr}

\author{Georg Menz}
\address{Department of Mathematics, University of California, Los Angeles}
\email{gmenz@math.ucla.edu}


\date{\today}

\begin{abstract}
We consider a one-dimensional lattice system of unbounded, real-valued spins with arbitrary strong, quadratic, finite-range interaction. We show the equivalence of the grand canonical ensemble (gce) and the canonical ensemble (ce), in the sense of observables and correlations. A direct consequence is that the correlations of the ce decay exponentially plus a volume correction term. The volume correction term is uniform in the external field, the mean spin and scales optimally in the system size. This extends prior results of Cancrini \& Martinelli for bounded discrete spins to unbounded continuous spins. The result is obtained by adapting Cancrini \& Martinelli's method combined with authors' recent approach on continuous real-valued spin systems.
\end{abstract}

\maketitle

\section{Introduction}

We are interested in studying the equivalence of the grand canonical ensemble (gce) and the canonical ensemble (ce). We consider a one-dimensional lattice system of unbounded real-valued spins denoted by~$\Lambda \subset \mathbb{Z}$. Throughout this paper we assume~$\Lambda$ is a finite system given by~$\{1, \cdots, N\}$ for convenience. In the Ising model the spin~$x_i$ at each site~$i \in \Lambda$ can take on the value 0 or 1. In this paper, we assume that the spins are real-valued and unbounded, i.e.,~$x_i \in \mathbb{R}$. A configuration of the lattice system is given by a vector~$x \in \mathbb{R}^{N}$ and the energy of a configuration~$x$ is given by the Hamiltonian~$H : \mathbb{R}^{N} \to \mathbb{R}$. For the detailed definition of the Hamiltonian~$H$ we refer to Section~\ref{s_setting_and_main_results}. \\

We consider two ensembles of the lattice system. The first ensemble is the grand-canonical ensemble (gce). For each~$\sigma \in \mathbb{R}$ we define the gce~$\mu^{\sigma}$ by the finite-volume Gibbs measure
\begin{align}
  \mu^\sigma(dx) = \frac{1}{Z} \exp\left( \sigma \sum_{i=1}^N x_i- H(x) \right) dx.
\end{align}
Here,~$Z$ is a generic normalization constant making the measure~$\mu^\sigma$ a probability measure. More precisely,~$Z$ is given by
\begin{align}
Z = \int_{\mathbb{R}^N} \exp \left( \sigma \sum_{i=1}^{N} x_i -H(x) \right) dx.
\end{align}
The constant~$\sigma \in \mathbb{R}$ is interpreted as an external field. The second ensemble is the canonical ensemble (ce), which emerges from the gce by conditioning on the mean spins. More precisely, for given fixed mean spin~$m \in \mathbb{R}$, the ce is given by the probability measure
\begin{align}
  \mu_m(dx)=  \frac{1}{Z} \  \mathds{1}_{\left\{ \frac{1}{N}\sum_{i=1}^N x_i= m \right\}  } (x)\ \exp\left(\sigma \sum_{i=1}^{N} x_i - H(x)\right) \mathcal{L}^{N-1}(dx),
\end{align} 
where~$\mathcal{L}^{N-1}$ denotes the~$(N-1)$-dimensional Hausdorff measure. In Section~\ref{s_setting_and_main_results}, we will revisit the definitions of gce and ce and see how external field~$\sigma$ and mean spin~$m$ are related. \\

There are many levels of defining equivalence of ensembles. In this article we follow the presentation of~\cite{MR00}, where three levels of equivalence of ensembles were introduced: on the level of thermodynamic functions, on the level of observables, and on the level of correlations. That is, in the thermodynamic limit (size~$N$ of the system goes to~$\infty$), free energy, expectation of intensive observable, and correlation of two intensive functions are independent of the ensemble used. In this article, we study equivalence of the grand canonical ensemble and the canonical ensemble on the level of observables and on the level of correlations. We refer to Section~\ref{s_setting_and_main_results} for more details. For further background, we refer the reader to~\cite{StZe91,LePfSu94,Ge95,Ad06,Tou15}. \\

The equivalence of ensembles in one-dimensional lattice system was well studied for discrete spin values (see~\cite{DoTi77}) or for quadratic Hamiltonians (see~\cite{Ge95}). Cancrini \& Martinelli~\cite{MR00} provided quantitative optimal upper bounds of the equivalence of the gce and the ce in the case of bounded discrete spins. The case where the spin values are unbounded real-valued and the Hamiltonian is not quadratic was studied by the authors (see~\cite{kwme18a}). However, although the results in~\cite{kwme18a} are quantitative the bounds are not optimal in terms of the system size. The question whether quantitative optimal bounds can be obtained in our setting remained open. In this article, we show that this is indeed the case. The main results of this article, i.e., Theorem~\ref{p_equivalence_observables} and Theorem~\ref{p_equivalence_decay_of_correlations}, state that the gce and ce are equivalent. The upper bounds are explicit and are scaling optimally. We therefore extend the results of~\cite{MR00} from bounded discrete spins to unbounded continuous spins, and at the same time improve the estimates of~\cite{kwme18a}. \\

In the proof of Theorem~\ref{p_equivalence_observables} we follow a Cancrini \& Martinelli's method combined with authors' recent approach on continuous real-valued spin systems (see for example~\cite{kwme18a}). Like in Cancrini \& Martinelli's method, we use Fourier transform to write
\begin{align} 
&\mathbb{E}_{\mu_m} \left[ f (X) \right] - \mathbb{E}_{\mu^\sigma} \left[ f (X) \right]\\
& = \frac{\int_{\mathbb{R}}\mathbb{E}_{\mu^{\sigma}} \left[ \left( f (X) - \mathbb{E}_{\mu^{\sigma}} \left[ f (X) \right]\right) \exp \left( i \frac{1}{\sqrt{N}} \sum_{i=1}^{N} \left(X_i -m_i \right) \xi\right) \right]d\xi}{\int_{\mathbb{R}}\mathbb{E}_{\mu^{\sigma}} \left[  \exp \left( i \frac{1}{\sqrt{N}} \sum_{i=1}^{N} \left(X_i -m_i \right) \xi\right) \right]d\xi}.
\end{align}
Then we first prove the theorem for an intensive function~$f$ which is "almost orthogonal" to the random variable~$\sum_{i \in \supp f} X_i$. That is, the covariance between~$f$ and~$\sum_{i \in \supp f }X_i$ is of order~$\frac{1}{N}$. For the general case, we decompose the intensive function~$f$ into "almost orthogonal" part and the remainder. We note that this can be done by subtracting~$C \sum_{i \in \supp f} X_i$ from~$f$, where~$C$ is a suitable constant depending on~$f$. Then by further decomposing the remainder~$C \sum_{i \in \supp f} X_i$ into "almost orthogonal remainder" and the rest, we obtain the desired result with a help of the identity
\begin{align}
\sum_{i =1}^{N} \mathbb{E}_{\mu^{\sigma}} \left[ X_i \right] = \sum_{i=1}^{N} \mathbb{E}_{\mu_m} \left[ X_i \right],
\end{align}
which will be verified in Section~\ref{s_setting_and_main_results}. The technical estimates follow and improve arguments of~\cite{kwme18a}. \\

The proof of Theorem~\ref{p_equivalence_decay_of_correlations} is also inspired by arguments of~\cite{kwme18a}. There, a similar result was deduced but the scaling is sub-optimal compared to the results of discrete case (cf.~\cite{MR00}). The main idea of the proof is to write the difference of correlations between two ensemble in terms of the expectations with respect to the gce. While the result of~\cite{kwme18a} rely on the second order Taylor expansion with third order moment bounds (cf. Lemma~\ref{l_third_moment}), which resulted in sub-optimality of the result, we make use of the third order Taylor expansion and the fourth order moment bounds (cf. Lemma~\ref{l_fourth_moment}) to achieve an optimal scaling. \\

The main results of this article, see Theorem~\ref{p_equivalence_observables} and Theorem~\ref{p_equivalence_decay_of_correlations}, also complement the recent results of Cancrini \& Olla~\cite{CaOl17}. In~\cite{CaOl17}, the equivalence of ensembles for \emph{extensive} observables was deduced in particle systems via an Edgeworth expansion, whereas our result applies to \emph{intensive} observables. In particular, the result of~\cite{CaOl17} is optimal proving the Lebowitz-Percus-Verlet formula and applies to a wider class of models. However, it is conditional, i.e., the assumptions are needed to be verified for any particular choice of observable~$f$, while our result is unconditional. It might be possible that their method could be extended to intensive observables but to get an unconditional result, one would need to prove their assumptions for a class of intensive variables. This is equivalent to proving a weaker version of equivalence of ensembles. It is not clear if this is possible. It also might be that for certain intensive observables in non-Gaussian models those assumptions fail. Our results, Theorem~\ref{p_equivalence_observables} and Theorem~\ref{p_equivalence_decay_of_correlations}, are unconditional and apply to wide class of intensive observables. \\

An important implication of the equivalence of ensemble is the decay of correlations of the ce (cf. Theorem~\ref{p_exponential_decay_of_correlations_of_canonical_ensemble}). In~\cite{kwme18a}, the decay of correlation was deduced under the same settings of this article. Because the equivalence of ensembles result in~\cite{kwme18a} was sub-optimal, also the decay of correlation result was sub-optimal. We revisit this statement with optimal scaling as a corollary of our main result, Theorem~\ref{p_equivalence_decay_of_correlations}. While the decay of correlation of an ensemble itself is a very interesting property, it also plays an integral role in deducing a uniform log-Sobolev inequality (LSI) of the ensemble. Indeed, in the case of strong, finite-range interactions, the uniform LSI of the ce was deduced in~\cite{KwMe18b}. Moreover, it was also shown with the help of decay of correlation that the ce on the one-dimensional lattice does not have a phase transition (see~\cite{kwme18a}).\\

Let us mention some open questions and problems:
\begin{itemize}
\item Instead of considering finite-range interaction, is it possible to deduce similar results for infinite-range, algebraically decaying interactions? More precisely, is it possible to extend the results of~\cite{MeNi14} from the gce to the ce? Is the same algebraic order of decay sufficient, i.e.~of the order~$2+ \varepsilon$, or does one need a higher order of decay? \\[-1ex]

\item Is it possible to consider more general Hamiltonians? For example, our argument is based on the fact that the single-site potentials are perturbed quadratic, especially when we use the results of~\cite{KwMe18}. One would like to have general super-quadratic potentials as was for example used in~\cite{MeOt13}. Also, it would be nice to consider more general interactions than quadratic or pairwise interaction.\\[-1ex]
	
\item Is it possible to generalize the results to vector-valued spin systems?.
\end{itemize}

We conclude the introduction by giving an overview over the article. In Section~\ref{s_setting_and_main_results}, we introduce the precise setting and present the main results. In Section~\ref{s_auxiliary_lemmas}, we provide several auxiliary results. In Section~\ref{s_proof_equiv_obs}, we prove the main results, i.e., the equivalence of the the gce and the ce on the level of observables and correlations (cf.~Theorem~\ref{p_equivalence_observables} and Theorem~\ref{p_equivalence_decay_of_correlations}). In Section~\ref{s_proof_1st_der_main_proposition} and Section~\ref{s_proof_2nd_der_main_proposition}, we provide the proof of additional ingredients that are needed in the proof of Theorem~\ref{p_equivalence_observables} and Theorem~\ref{p_equivalence_decay_of_correlations}, respectively. \\

\section*{Conventions and Notation}

\begin{itemize}
\item The symbol~$T_{(k)}$ denotes the term that is given by the line~$(k)$.
\item We denote with~$0<C<\infty$ a generic uniform constant. This means that the actual value of~$C$ might change from line to line or even within a line.
\item Uniform means that a statement holds uniformly in the system size~$N$, the mean spin~$m$ and the external field~$s$.
\item $a \lesssim b$ denotes that there is a uniform constant~$C$ such that~$a \leq C b$.
\item $a \sim b$ means that~$a \lesssim b$ and~$b \lesssim a$.
\item $\mathcal{L}^{k}$ denotes the $k$-dimensional Hausdorff measure. If there is no cause of confusion we write~$\mathcal{L}$.
\item $Z$ is a generic normalization constant. It denotes the partition function of a measure.  
\item For each~$N \in \mathbb{N}$,~$[N]$ denotes the set~$\left\{ 1, \ldots N \right\}$.
\item For a vector~$x \in \mathbb{R}^{[N]}$ and a set~$A \subset [N]$,~$x^A \in \mathbb{R}^{A}$ denotes the vector $ (x^A)_{i} = x_i$ for all~$i \in A$.
\item For a function~$f : \mathbb{R}^{[N]} \to \mathbb{C}$, denote~$\supp f $ by the minimal subset of~$\mathbb{Z}$ with $f(x) = f\left(x^{\supp f} \right)$. 
\item For a vector~$x \in \mathbb{R}^{n}$,~$| x|$ denotes the standard Euclidean norm of~$x$.
\item For a function~$f : \mathbb{R}^{[N]}$, the~$L^p$ norm of~$f$ with respect to gce~$\mu^{\sigma}$ is given by
\begin{align}
    \| \nabla f \|_{L^p (\mu^{\sigma})} = \left(  \mathbb{E}_{\mu^{\sigma}} \left[ | \nabla f |^p \right]  \right)^ \frac{1}{p}.
\end{align}
$L^{\infty}$ norm is~$\| \nabla f \|_{\infty}:= \| \nabla f \|_{L^{\infty} (\mu^{\sigma})} : = | g | _{L^{\infty} (\mu^{\sigma})}$, where~$ g = |\nabla f |$.
\end{itemize}

\medskip

\section{Setting and main results}
\label{s_setting_and_main_results}

We consider a system of unbounded continuous spins on the sublattice~$\{1, \cdots, N \} \subset \mathbb{Z}$. The Hamiltonian $H:\mathbb{R}^{N} \to \mathbb{R}$ of the system is defined as
\begin{align}\label{e_d_hamiltonian}
H(x) &= \sum_{i = 1 }^{N} \left( \psi (x_i) + s_i x_i +\frac{1}{2}\sum_{j : \ 1 \leq |j-i| \leq R } M_{ij}x_i x_j \right) \\
& = \sum_{i =1} ^{N} \left( \psi_b (x_i) + s_ix_i + \frac{1}{2} \sum_{j : |j-i| \leq R} M_{ij}x_i x_j  \right),
\end{align}
where~$\psi(z) : = \frac{1}{2}z^2 + \psi_b (z)$ and~$M_{ii} : = 1$.
We assume the following:
\begin{itemize}
\item The function~$\psi_b: \mathbb{R} \to \mathbb{R}$ satisfies
 \begin{align}\label{e_nonconvexity_bounds_on_perturbation}
 |\psi_b|_{\infty} + |\psi'_b|_{\infty}  + |\psi''_b|_{\infty} < \infty. 
 \end{align}
It is best to imagine~$\psi(z) = \frac{1}{2}z^2 + \psi_b(z)$ as a double-well potential (see Figure~\ref{f_double_well}).  
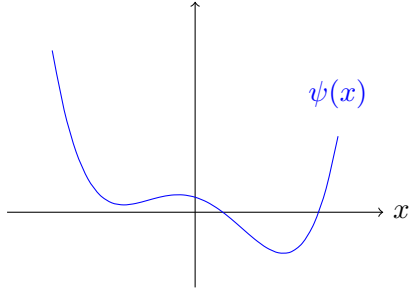
\begin{figure}[t]
\centering
\begin{tikzpicture}
      \draw[->] (-2.5,0) -- (2.5,0) node[right] {$x$};
      \draw[->] (0,-1) -- (0,2.8) node[above] {};
      \draw[scale=1,domain=-1.9:1.9,smooth,variable=\x,blue] plot ({\x},{.3*\x*\x*\x*\x-.7*\x*\x-.3*\x+.2});
      
\node[align=center, below, blue] at (1.9, 1.9) {$\psi(x)$};      
\end{tikzpicture}
\caption{Example of a single-site potential~$\psi$}\label{f_double_well}
\end{figure}

\item The interaction is symmetric i.e. 
\begin{align}
M_{ij} = M_{ji}  \qquad &\text{for all distinct} \ i, j \in \mathbb{Z}. 
\end{align}

\item The fixed, finite number~$R \in \mathbb{N}$ models the range of interactions between the particles in the system i.e.~it holds that~$M_{ij}=0$ for all~$i,j$ such that~$|i-j|> R$ or~$(i, j) \notin \{ 1, \cdots, N \} \times \{ 1, \cdots, N \}$. \\

\item The matrix~$M=(M_{ij})$ is strictly diagonal dominant i.e. for some~$\delta>0$, it holds for any~$i \in \mathbb{Z}$ that
\begin{align}
\sum_{ 1 \leq |j-i| \leq R} |M_{ij}| + \delta \leq M_{ii} = 1. 
\end{align}

\item The vector~$s = (s_i) \in \mathbb{R}^N$ is arbitrary. It models the interaction with an inhomogeneous external field. Because the interaction is quadratic, this term also models the interaction of the system with the boundary. 
\end{itemize}

\medskip 
\begin{definition} [The canonical and grand canonical ensemble] The gce~$\mu^{\sigma}$ associated to the Hamiltonian~$H$ is the probability measure on~$\mathbb{R}^{N}$ given by the Lebesgue density
\begin{align} \label{d_gc_ensemble}
\mu^{ \sigma} \left(dx\right) : = \frac{1}{Z} \exp\left( \sigma \sum_{ i =1}^{N} x_i  - H(x ) \right) dx, 
\end{align}
where~$dx$ denotes the Lebesgue measure on~$\mathbb{R}^N$. The ce~$\mu_m $ is the probability measure on
\begin{align}\label{e_d_X_lambda_m}
X_{N, m} : = \left \{ x \in \mathbb{R}^N : \ \frac{1}{N} \sum_{i =1}^N x_i =m  \right \} \subset \mathbb{R}^{N}
\end{align}
with density
\begin{align} \label{d_ce}
\mu_m  (dx) : = \frac{1}{Z} \mathds{1}_{ \left\{ \frac{1}{N} \sum_{i =1}^{N} x_i =m \right\}}\left(x\right) \exp\left( - H(x) \right) \mathcal{L}^{N-1}(dx), \label{d_ce}
\end{align}
where~$\mathcal{L}^{N-1}(dx)$ denotes the~$(N-1)$-dimensional Hausdorff measure supported on~$X_{N, m}$. 
\end{definition}

\medskip 

\begin{remark} The ce~$\mu_m $ emerges from the gce~$\mu^{\sigma}$ by conditioning on the mean spin
\begin{align} \label{e_spin_restriction}
\frac{1}{N} \sum_{ i =1}^{N} x_i = m.
\end{align}
More precisely, given~\eqref{e_spin_restriction}, the term~$\sigma \sum_{i=1}^{N} x_i$ inside the exponential in~\eqref{d_gc_ensemble} acts like a constant and hence is cancelled out with the normalization constant~$Z$ as follows:
\begin{align}
\mu^{\sigma}  \left(dx   \mid \frac{1}{N}\sum_{i =1}^{N} x_i =m \right) & = \frac{1}{Z} \mathds{1}_{ \left\{ \frac{1}{N} \sum_{i=1}^{N} x_i =m \right\}}(x) \exp\left(\sigma m N - H(x) \right) \mathcal{L}^{N-1}(dx)  \\
& = \frac{1}{\widetilde{Z}} \mathds{1}_{ \left\{ \frac{1}{N} \sum_{i =1}^{N} x_i =m \right\}}\left(x\right) \exp\left( - H(x) \right) \mathcal{L}^{N-1}(dx) \\
& = \mu_m (dx).
\end{align}
We note that the ce~$\mu_m $ does not have a dependence on~$\sigma$ anymore, even though it emerged from the gce~$\mu^{\sigma}$.
\end{remark}

\medskip 

To relate the external field~$\sigma$ of~$\mu^{\sigma}$ and the mean spin~$m$ of~$\mu_m $ we further introduce following definition. \\

\begin{definition}[The free energy of the gce] The free energy~$A_{gce} : \mathbb{R} \to \mathbb{R}$ of the gce~$\mu^{ \sigma}$ is defined as
\begin{align}
A_{gce} (\sigma) : = \frac{1}{N}\ln \int_{\mathbb{R}^{N}} \exp \left( \sigma \sum_{i=1}^{N} x_i - H(x) \right)dx. \label{e_def_a_gce}
\end{align}
\end{definition}

\medskip 

Let~$X= (X_i)_{i=1}^{N}$ be a random variable distributed according to the gce~$\mu^{\sigma}$. A direct calculation yields
\begin{align}
\frac{d^2}{d\sigma^2} A_{gce} (\sigma) = \frac{1}{N} \var_{\mu^{ \sigma}} \left(  \sum_{i=1}^{N} X_i \right).
\end{align}

\medskip 
Later on, we will use the following observation. \\

\begin{lemma}[Lemma 3.1 in~\cite{KwMe18b}] \label{l_variance_estimate} There is a constant~$C \in (0, \infty)$ which is uniform in the system size~$N$ and the external fields~$s$,~$\sigma$ such that
\begin{align}
\frac{1}{C} \leq \frac{1}{N} \var_{\mu^{ \sigma}} \left( \sum_{i=1}^{N} X_i  \right)\leq C.
\end{align}
\end{lemma}
 
\medskip 

\begin{corollary} The free energy~$A_{gce}$ of the gce~$\mu^{\sigma}$ is strictly convex in the sense that there is a constant~$C \in (0, \infty)$ independent of the system size~$N$ and the external fields~$s$,~$\sigma$ such that
\begin{align}
\frac{1}{C} \leq \frac{d^2}{d\sigma^2} A_{gce} \leq C.
\end{align}
\end{corollary}
\medskip 

We denote with~$\mathcal{H}_N$ the Legendre transform of the free energy~$A_{gce}$. It is defined by
\begin{align}
\mathcal{H}_N (m) = \sup_{\sigma \in \mathbb{R}} \left( \sigma m - A_{gce}(\sigma) \right).
\end{align}
Due to the strict convexity of the free energy~$A_{gce}$, it follows that for each~$m \in \mathbb{R}$ there is a unique~$\sigma = \sigma(m) \in \mathbb{R}$ such that
\begin{align} \label{e_relation_sigma_m}
\mathcal{H}_N (m) = \sigma(m) m - A_{gce} (\sigma(m)).
\end{align}
From now on, we always assume that~$\sigma$ and~$m$ is related by~\eqref{e_relation_sigma_m}. In particular, it holds that
\begin{align}
\frac{d}{d\sigma} A_{gce} (\sigma) = m
\end{align}
and
\begin{align}
\frac{d}{dm} \mathcal{H}_N (m) = \frac{d \sigma}{dm} m + \sigma(m) - \frac{d}{d\sigma} A_{gce} (\sigma(m)) \frac{d\sigma}{dm} = \sigma(m).
\end{align}
Setting~$m_i := \int x_i \mu^{ \sigma}\left(dx\right)$ yields
\begin{align}
m &= \frac{d}{d\sigma} A_{gce} (\sigma) \\
& = \frac{1}{N} \frac{ \int_{\mathbb{R}^N} \sum_{i=1}^{N} x_i \exp \left( \sigma \sum_{i=1}^{N} x_i - H(x) \right)dx }{  \int_{\mathbb{R}^N}  \exp \left( \sigma \sum_{i=1}^{N} x_i - H(x) \right) dx } \\
& = \frac{1}{N} \sum_{i=1}^{N} m_i. \label{e_m=sum_mi}
\end{align}



\medskip  





\begin{definition} [Local, intensive, and extensive functions/ observables]
For a function~$f : \mathbb{R}^{\mathbb{Z}} \to \mathbb{C}$, denote~$\supp f $ by the minimal subset of~$\mathbb{Z}$ with $f(x) = f\left(x^{\supp f} \right)$. We call~$f$ a local function if it has a finite support independent of~$N$. A function~$f$ is called intensive if there is a positive constant~$\varepsilon$ such that $|\supp f | \lesssim N^{1-\varepsilon}$. A function~$f$ is called extensive if it is not intensive.
\end{definition}

\medskip 

Let us now turn to the first main result of this article, the equivalence of the ce and gce on the level of observables.  \\

\begin{theorem}[Equivalence of the ce and gce on the level of observables]  \label{p_equivalence_observables}
Let~$f: \mathbb{R}^{\mathbb{Z}} \to \mathbb{R}$ be an intensive function. There are constants~$C \in (0, \infty)$ and~$N_0 \in \mathbb{N}$ independent of the external field~$s$ and the mean spin~$m$ such that for all~$N \geq N_0$, it holds that
\begin{align}
\left| \mathbb{E}_{\mu^{ \sigma}}\left[ f\right] - \mathbb{E}_{\mu_m} \left[ f \right] \right| \leq C \frac{|\supp f| }{N}\| \nabla f \|_{\infty}. \label{e_equivalence_observables}
\end{align}
\end{theorem}

\medskip 

We provide the proof of Theorem~\ref{p_equivalence_observables} in Section~\ref{s_proof_equiv_obs}. \\

Let us now turn to the second main result of this article: \\

\begin{theorem}[Equivalence of the ce and gce on the level of correlations] \label{p_equivalence_decay_of_correlations}
Let~$f, g : \mathbb{R}^{N} \to \mathbb{R}$ be intensive functions. There exist constants~$C \in (0, \infty)$ and~$N_0 \in \mathbb{N}$ independent of the external field~$s$ and the mean spin~$m$ such that for all~$N \geq N_0$, it holds that
\begin{align}
& \left|  \cov_{\mu_m}  (f ,  g) - \cov_{\mu^{\sigma}}  (f  ,  g) \right| \\
&\qquad \leq C \  \| \nabla f \|_{\infty}\| \nabla g \|_{\infty} \left( \frac{  |\supp f | + |\supp g | }{N} + \exp\left(-C\text{dist}\left( \supp f, \supp g \right) \right) \right)  . \label{e_equiv_dec_cor}
\end{align}
\end{theorem}

\medskip 

We give the proof of Theorem~\ref{p_equivalence_decay_of_correlations} in Section~\ref{s_proof_equiv_obs}. \\

A consequence of Theorem~\ref{p_equivalence_decay_of_correlations} is the decay of correlations of the ce. For that purpose let us recall that for one-dimensional lattice systems the correlations of the gce decay exponentially fast (\cite[Lemma 6]{KwMe18}. See also~\cite[Theorem 1.4]{MeNi14}). \\

\begin{theorem}[Lemma 6 in~\cite{KwMe18}] \label{p_decay_of_correlations_gce}
Let~$f, g : \mathbb{R}^{N} \to \mathbb{R}$ be intensive functions. Then 
\begin{align}
&\left| \cov_{\mu^{\sigma}} \left( f, g \right) \right| \lesssim \|\nabla f\|_{L^2 (\mu^{\sigma})}\|\nabla g\|_{L^2 (\mu^{ \sigma})} \exp \left( -C \text{dist} \left( \supp f , \supp g \right) \right). \label{e_decay_of_correlation_ce}
\end{align}
\end{theorem}

\medskip 

A straightforward combination of Theorem~\ref{p_decay_of_correlations_gce} and Theorem~\ref{p_equivalence_decay_of_correlations} yields another main result of this article: \\

\begin{theorem}[Decay of correlations of the ce]\label{p_exponential_decay_of_correlations_of_canonical_ensemble}
Under the same assumptions as in Theorem~\ref{p_equivalence_decay_of_correlations}, it holds that
\begin{align} \label{e_decay_correlations_ce}
\left|\cov_{\mu_m } \left( f, g \right) \right| \leq C \ \| \nabla f \|_{\infty}\| \nabla g \|_{\infty}  \left( \frac{ |\supp f | + |\supp g | }{N} + \exp\left(-C\text{dist}\left( \supp f, \supp g \right) \right) \right).
\end{align}
\end{theorem}
\medskip 

\begin{remark}
One should compare Theorem~\ref{p_equivalence_observables} and Theorem~\ref{p_equivalence_decay_of_correlations} with~\cite[Theorem 2]{kwme18a} and~\cite[Theorem 5]{kwme18a}, respectively. There, similar results were deduced under the same settings. The scaling of system size~$N$ are improved from~$N^{ \frac{1}{2} - \varepsilon}$ and~$N^{1 - \varepsilon}$ to~$N$, where~$\varepsilon$ is a positive constant. This is consistent with the result presented in~\cite{MR00}, where equivalence of observables and correlations were deduced in the discrete spin system (see Theorem 4.1 and Proposition 7.3 in~\cite{MR00}). The main difference is that we use~$L^{\infty}$ norm of~$\nabla f$, while~\cite{MR00} used~$L^{\infty}$ norm of~$f$. However, such bounds with~$\| f \|_{\infty}$ have limited use in the continuous settings. For example, one could not deduce the decay of the spin-spin correlation function (see Theorem~\ref{p_decay_two_point_ce}).
\end{remark}

\medskip 
Let us now illustrate the use of Theorem~\ref{p_exponential_decay_of_correlations_of_canonical_ensemble} by deducing the decay of the spin-spin correlation function of the canonical ensemble. \\

\begin{corollary}[Decay of the spin-spin correlation function of the ce]\label{p_decay_two_point_ce}
There exist constants~$N_0 \in \mathbb{N}$ and~$C \in (0, \infty)$ independent of the external field~$s$ and the mean spin~$m$ such that for any~$N \geq N_0$, it holds that for any~$i, j \in \{ 1, \cdots, N\}$,
\begin{align}
\left| \cov_{\mu_m} \left( X_i, X_j \right) \right| \leq C \left( \frac{1}{N} + \exp \left( - C |i-j| \right)\right). \label{e_can_exp_dec}
\end{align}
\end{corollary}

\medskip 

\begin{remark}
Compared to Theorem~\ref{p_decay_of_correlations_gce}, there appears an additional volume correction term~$\frac{1}{N}$ in Corollary~\ref{p_decay_two_point_ce}. This term is due to the mean constraint~$\frac{1}{N} \sum_{i=1}^{N} x_i = m$ and is optimal. For example, assuming that the Hamiltonian~$H$ is symmetric, we have
\begin{align}
\cov_{\mu_m} (X_1, X_2 ) = \cov_{\mu_m} (X_i , X_j) \qquad \text{ for all distinct } i, j \in \{1, \cdots, N\}.    
\end{align}
Thus we get
\begin{align}
\cov_{\mu_m} (X_1, X_2) &= \frac{1}{N-1} \cov_{\mu_m} (X_1, X_2 + \cdots X_N )\\
& = \frac{1}{N-1} \cov_{\mu_m} (X_1, Nm - X_1)\\
& = - \frac{1}{N-1} \var_{\mu_m} (X_1).
\end{align}
\end{remark}

\medskip 

For the proofs of the main results of this article, we refer to Section~\ref{s_proof_equiv_obs}.

\section{Auxiliary Lemmas} \label{s_auxiliary_lemmas}

\subsection{Basic Properties of the gce~$\mu^{\sigma}$ and ce~$\mu_m$}

In this section we provide auxiliary estimates that will be needed in the proof of Theorem~\ref{p_equivalence_observables} and Theorem~\ref{p_equivalence_decay_of_correlations}. \\

Let~$g$ be the density of the random variable
\begin{align}
\frac{1}{\sqrt{N}} \sum_{i =1}^{N} \left( X_i - m \right) \overset{\eqref{e_m=sum_mi}}{=} \frac{1}{\sqrt{N}} \sum_{i=1}^{N} \left( X_i - m_i \right),
\end{align}
where the random vector~$X=(X_i)_{i=1}^{N}$ is distributed according to~$\mu^{\sigma}$.
The following proposition provides estimates for~$g(0)$. \\

\begin{proposition} [Proposition 1 in~\cite{KwMe18}] \label{p_main_computation}
For each~$\alpha>0$ and~$\beta > \frac{1}{2}$, there exist constants~$C \in (0, \infty)$ and~$N_0 \in \mathbb{N}$ independent of the external field~$s$ and the mean spin~$m$ such that for all~$N \geq N_0$, it holds that
\begin{align}
\frac{1}{C} \leq g(0) \leq C , \qquad 
\left| \frac{d}{d\sigma} g(0) \right| \lesssim N^{\alpha}  \qquad \text{and} \qquad 
\left| \frac{d^2}{d\sigma^2} g(0) \right| \lesssim N^{\beta}.
\end{align}
Moreover, it holds that
\begin{align}
\frac{1}{C} \leq \int_{\mathbb{R}} \mathbb{E}_{\mu^{\sigma}} \left[ \exp \left( \frac{1}{\sqrt{N}} \sum_{i=1}^{N} \left( X_i -m_i \right) \xi \right)  \right] d \xi \leq C.
\end{align}
\end{proposition}
\medskip

In the remaining section we provide several estimates for moments of observables and correlations. Those estimates are a bit exhaustive but a standard ingredient when studying equivalence of ensembles. The following lemma provides a general moment estimate for the gce. \\

\begin{lemma} \label{l_function_moment_estimate}
For each~$k \geq 1$, there is a constant~$C = C(k)$ such that for any smooth function~$f : \mathbb{R}^{\Lambda} \to \mathbb{R}$
\begin{align} \label{e_function_moment_estimate}
\mathbb{E}_{\mu^{ \sigma}} \left[ \left|f(X) - \mathbb{E}_{\mu^{ \sigma}} \left[f(X) \right] \right|^k \right] \leq C(k) \| \nabla f \|_{\infty}^k.
\end{align}
\end{lemma}
\medskip 

The statement of Lemma~\ref{l_function_moment_estimate} is a simple extension of (23) in~\cite{KwMe18} from the special case of~$f(x) = x_i$ to general functions. For the convenience of the reader, we restate the short argument. \\

\noindent \emph{Proof of Lemma~\ref{l_function_moment_estimate}.} \ It is well known that the gce~$\mu^{\sigma}$ satisfies a uniform LSI and Poincar\'e inequality (see for example~\cite{HeMe16}). The case~$k=2$ easily follows from an application of Poincar\'e inequality. More precisely, we have
\begin{align}
\mathbb{E}_{\mu^{ \sigma}} \left[ \left|f(X) - \mathbb{E}_{\mu^{ \sigma}} \left[f(X) \right] \right|^2 \right] \leq \frac{1}{\rho} \mathbb{E}_{\mu^{\sigma}} \left[ |\nabla f |^2  \right]  \leq \frac{1}{\rho} \| \nabla f \|_{\infty} ^2,
\end{align}
where~$\rho>0$ is a uniform constant in Poincar\'e inequality. Thanks to the Schwarz inequality,~\eqref{e_function_moment_estimate} also holds for~$k=1$. Assume that~\eqref{e_function_moment_estimate} holds for some~$k=2n \geq 2$. Again, Poincar\'e inequality implies that
\begin{align}
&\mathbb{E}_{\mu^{ \sigma}} \left[ \left|f(X) - \mathbb{E}_{\mu^{ \sigma}} \left[f(X) \right] \right|^{2n+2} \right] - \left(\mathbb{E}_{\mu^{ \sigma}} \left[ \left|f(X) - \mathbb{E}_{\mu^{ \sigma}} \left[f(X) \right] \right|^{n+1} \right] \right)^2 \\
&\qquad  \leq \frac{1}{\rho} \mathbb{E}_{\mu^{\sigma}} \left[ \left| \nabla \left( \left|f(X) - \mathbb{E}_{\mu^{ \sigma}} \left[f(X) \right] \right|^{n+1} \right) \right|^2 \right] \\
& \qquad \leq \frac{n+1}{\rho} \| \nabla f \|_{\infty}^2 \mathbb{E}_{\mu^{\sigma}}  \left[ \left| f(X) - \mathbb{E}_{\mu^{ \sigma}} \left[f(X) \right] \right|^{2n} \right] \\
& \qquad \leq \frac{n+1}{\rho} \| \nabla f \|_{\infty}^{2n+2}. \label{e_proof_function_moment1}
\end{align}
Because~$n+1 \leq 2n$, Schwarz inequality implies
\begin{align}
\left(\mathbb{E}_{\mu^{ \sigma}} \left[ \left|f(X) - \mathbb{E}_{\mu^{ \sigma}} \left[f(X) \right] \right|^{n+1} \right] \right)^2 &\leq  \mathbb{E}_{\mu^{\sigma}} \left[ \left|f(X) - \mathbb{E}_{\mu^{ \sigma}} \left[f(X) \right] \right|^{2n} \right]^{\frac{n+1}{2n} \cdot 2 } \\
& \leq \left(  C(2n) \| \nabla f \|_{\infty}^{2n} \right)^{\frac{n+1}{n}} \lesssim \| \nabla f \|_{\infty}^{2n+1}. \label{e_proof_function_moment2}
\end{align}
A combination of~\eqref{e_proof_function_moment1} and~\eqref{e_proof_function_moment2} proves~\eqref{e_function_moment_estimate} for~$k = 2n+2$. Schwarz inequality also implies that this holds for~$k = 2n+1$. Then mathematical induction concludes the proof of Lemma~\ref{l_function_moment_estimate}. \qed 

\medskip

The next statement is a direct consequence of Lemma~\ref{l_function_moment_estimate}. \\

\begin{corollary} \label{l_moment_estimate} For each~$i \in [N]$, we define
\begin{align} \label{d_def_mi}
m_i : = \int x_i \mu^{ \sigma} \left(dx\right).
\end{align}
Then for each $k \geq 1$, there is a constant~$C=C(k)$ such that for each~$i \in [N]$
\begin{align}
\mathbb{E}_{\mu^{\sigma}}\left[ \left| X_i -m_i \right|^k \right] \leq C(k) .
\end{align}
\end{corollary}
 
\medskip

The next statement is an estimation of cubic moments. \\

\begin{lemma} \label{l_third_moment} Let~$(X_1, \cdots, X_n)$ be a real-valued random variable distributed according to the gce~$\mu^{\sigma}$. For each~$i \in \{1, \cdots, N\}$ denote~$Y_i := X_i -\mathbb{E}_{\mu^{\sigma}} \left[ X_i \right]$. Then for any subset~$A$ of $[N]$, it holds that
\begin{align}
\left| \mathbb{E}_{\mu^{\sigma}} \left[ \left( \sum_{i \in A} Y_i \right)^3 \right] \right| \lesssim |A|.
\end{align}
\end{lemma}

\medskip

\noindent \emph{Proof of Lemma~\ref{l_third_moment}.} \ We prove the case when~$A = \{1, \cdots, N\}$. General case follows from the same argument. For each pair~$(i, j , k) \subset \{1, \cdots, N\}$ with~$i \leq j \leq k$, we have by Theorem~\ref{p_decay_of_correlations_gce} and Corollary~\ref{l_moment_estimate} that
\begin{align}
\left| \mathbb{E}_{\mu^{\sigma}} \left[ Y_i Y_j Y_k \right] \right|  = \left| \cov_{\mu^{\sigma}} \left(  Y_i, Y_j Y_k \right)  \right| \lesssim \exp \left( - C |i-j| \right).
\end{align}
Similarly, one also gets
\begin{align}
\left| \mathbb{E}_{\mu^{\sigma}} \left[ Y_i Y_j Y_k \right] \right| \lesssim \exp \left( -C |j-k| \right)
\end{align}
and conclude
\begin{align} \label{e_third_key_estimate}
\left| \mathbb{E}_{\mu^{\sigma}} \left[ Y_i Y_j Y_k \right] \right| \lesssim \exp \left( -C \max \left( |i-j|, |j-k| \right)  \right).
\end{align}
Combined with the triangle inequality, the estimate~\eqref{e_third_key_estimate} yields
\begin{align}
\left| \mathbb{E}_{\mu^{\sigma}} \left[ \left( \sum_{i=1}^N Y_i \right)^3 \right] \right|  \lesssim \sum_{i \leq j \leq k} \mathbb{E} \left[ Y_i Y_j Y_k \right]  \lesssim \sum_{i \leq j \leq k } \exp \left( - C \max \left( |i-j|, |j-k| \right)  \right) .
\end{align}
For fixed~$j \in \{1, \cdots, N\}$ and~$d = \max \left( |i-j|, |j-k| \right)$, there are at most~$2 (d+1)$ pairs of~$(i,j,k)$ with~$i \leq j \leq k $. Therefore we conclude that
\begin{align}
\sum_{i \leq j \leq k} \exp \left( - C \max \left( |i-j|, |j-k| \right)\right) \leq \sum_{j=1}^{N} \sum_{d=0}^{N} 2(d+1) \exp \left( - Cd \right) \lesssim N.
\end{align}
This finishes the proof of Lemma~\ref{l_third_moment}.
\qed

\medskip

Next statement provides estimates for quartic moments. \\

\begin{lemma} \label{l_fourth_moment}
Under the same assumptions as in Lemma~\ref{l_third_moment}
\begin{align}
\mathbb{E}_{\mu^{\sigma}} \left[ \left( \sum_{i \in A} Y_i \right)^4 \right] \sim |A|^2.
\end{align}
\end{lemma}

\medskip

\noindent \emph{Proof of Lemma~\ref{l_fourth_moment}.} \ Again, we assume~$A= \{1, \cdots, N\}$. For each pair~$(i,j,k,l) \subset \{1, \cdots, N\}$ with~$i \leq j \leq k \leq l$ we have by Theorem~\ref{p_decay_of_correlations_gce} and Corollary~\ref{l_moment_estimate} that
\begin{align}
\left| \mathbb{E}_{\mu^{\sigma}} \left[ Y_i Y_j Y_k Y_l \right] \right| \lesssim \exp \left( -C \max \left( |i-j|, |k-l| \right) \right).
\end{align}
For fixed~$j \leq k$ and~$d = \max(|i-j|, |k-l|)$, there are at most~$2(d+1)$ pairs of~$(i, j, k, l)$ with~$i\leq j \leq k \leq l$. Thus we conclude
\begin{align}
\mathbb{E}_{\mu^{\sigma}} \left[ \left( \sum_{i=1}^N Y_i \right)^4 \right] \lesssim \sum_{i \leq j \leq k \leq l} \mathbb{E}_{\mu^{\sigma}} \left[ Y_i Y_j Y_k Y_l \right] \lesssim \sum_{j \leq k} \sum_{d=0}^{N-1} 2(d+1) \exp \left( -Cd \right) \lesssim N^2.
\end{align}

For the lower bound, we apply Lemma~\ref{l_variance_estimate} to get
\begin{align}
\mathbb{E}_{\mu^{\sigma}} \left[ \left( \sum_{i=1}^N Y_i \right)^4 \right] &= \var \left( \left( \sum_{i=1}^N Y_i \right)^2 \right) + \left( \mathbb{E}_{\mu^{\sigma}} \left[ \left( \sum_{i=1}^N Y_i \right)^2 \right] \right)^2 \\
& \geq \left( \mathbb{E}_{\mu^{\sigma}} \left[ \left( \sum_{i=1}^N Y_i \right)^2 \right] \right)^2 \gtrsim N^2.
\end{align}
\qed

\medskip

The next two statements provide estimates of covariance between an intensive function~$f$ and sum of spins. \\

\begin{lemma} \label{l_function_first_estimate}
Let~$f$ be an intensive function and~$A$ be any subset of~$[N]$. Then it holds that
\begin{align}
\left|\cov_{\mu^{\sigma}} \left( f(X) \ , \ \sum_{i \in A} X_i \right) \right|  \lesssim \| \nabla f \|_{L^2(\mu^{\sigma})} |\supp f |^{\frac{1}{2}}. \label{e_function_first_estimate}
\end{align}
\end{lemma}

\medskip 

\noindent \emph{Proof of Lemma~\ref{l_function_first_estimate}.} \ As before we assume~$A = [N]$. For each~$d \in \mathbb{N}$ denote~$S_d$ by
\begin{align} \label{d_def_sd}
S_d : = \{ k \in [N] : \text{dist}( \supp f \ , \ k ) = d \}.
\end{align}
We note that for each~$d \in \mathbb{N}$, the cardinality of~$S_d$ is bounded from above by~$2 |\supp f|$. By triangle inequality we have
\begin{align}
&\left|\cov_{\mu^{\sigma}} \left( f(X)  ,  \sum_{i=1}^N X_i \right) \right| \leq \left| \cov_{\mu^{\sigma}} \left( f(X)  ,  \sum_{k \in \supp f}  X_k \right)  \right|  + \sum_{d \geq 1} \left| \cov_{\mu^{\sigma}} \left( f(X)  ,  \sum_{k \in S_d}  X_k \right)  \right|. \label{e_function_first_estimate_decomposition}
\end{align}
Then a combination of Schwarz inequality, Poincar\'e inequality, Lemma~\ref{l_function_moment_estimate}, and Corollary~\ref{l_moment_estimate} yields that the first term in~\eqref{e_function_first_estimate_decomposition} is bounded by
\begin{align}
\left| \cov_{\mu^{\sigma}} \left( f(X)  ,  \sum_{k \in \supp f}  X_k \right)  \right| &\leq \| f(X) - \mathbb{E}_{\mu^{\sigma}}\left[ f(X) \right] \|_{L^2 (\mu^{\sigma})}  \left\| \sum_{k \in \supp f} \left( X_k -m_k \right)  \right\|_{L^2 (\mu^{\sigma})} \\
& \lesssim \| \nabla f \|_{L^2 (\mu^{\sigma})} |\supp f|^{\frac{1}{2}}.
\end{align}
Next, an application of Theorem~\ref{p_decay_of_correlations_gce} yields that the second term in~\eqref{e_function_first_estimate_decomposition} can be estimated as follows:
\begin{align}
\sum_{d \geq 1}  \left| \cov_{\mu^{\sigma}} \left( f(X)  ,   \sum_{k \in S_d} X_k \right)  \right| &\lesssim \sum_{d \geq 1} \| \nabla f \|_{L^2 (\mu^{\sigma})} |S_d|^{\frac{1}{2}} \exp \left( -Cd \right) \\
& \lesssim \| \nabla f \|_{L^2 (\mu^{\sigma})} |\supp f|^{\frac{1}{2}}\sum_{d \geq 1}  \exp \left( -Cd \right) \\
& \lesssim \| \nabla f \|_{L^2 (\mu^{\sigma})} |\supp f|^{\frac{1}{2}}.
\end{align}
\qed

\medskip 

\begin{lemma} \label{l_function_second_estimate}
Let~$f$ be an intensive function and~$A$ be any subset of~$[N]$. Then it holds that
\begin{align} \label{e_function_sec_est}
\left|\cov_{\mu^{\sigma}} \left(  f(X)  ,  \left( \sum_{i \in A} \left( X_i - m_i \right) \right)^2 \right)\right| \lesssim \| \nabla f \|_{L^4 (\mu^{\sigma} )} |\supp f|.
\end{align}
\end{lemma}

\medskip 

\noindent \emph{Proof of Lemma~\ref{l_function_second_estimate}.} \ Proof of Lemma~\ref{l_function_second_estimate} is motivated by the proof of Lemma~\ref{l_third_moment}. We use similar idea accompanied with more careful estimate when applying Theorem~\ref{p_decay_of_correlations_gce} and Corollary~\ref{l_moment_estimate}. In this proof, the set~$S$ denotes~$\supp f$. We first decompose the left hand side of~\eqref{e_function_sec_est} by
\begin{align}
&\cov_{\mu^{\sigma}} \left(  f(X)  ,  \left( \sum_{i \in A} \left( X_i - m_i \right) \right)^2 \right) \\
& \qquad  = \mathbb{E}_{\mu^{\sigma}} \left[ \left( f(X) - \mathbb{E}_{\mu^{\sigma}}\left[ f(X) \right] \right) \left(  \sum_{i \in S} \left( X_i -m_i \right)\right)^2 \right] \label{e_function_sec_1} \\
& \qquad \quad + 2 \mathbb{E}_{\mu^{\sigma}} \left[ \left( f(X) - \mathbb{E}_{\mu^{\sigma}}\left[ f(X) \right] \right) \sum_{i \in S} \left( X_i -m_i \right) \sum_{j \notin S} \left( X_j - m_j \right)   \right] \label{e_function_sec_2}\\
& \qquad \quad + \mathbb{E}_{\mu^{\sigma}} \left[ \left( f(X) - \mathbb{E}_{\mu^{\sigma}}\left[ f(X) \right] \right) \left(  \sum_{j \notin S} \left( X_i -m_i \right)\right)^2 \right]. \label{e_function_sec_3}
\end{align}
We estimate each term by term. Let us begin with estimating~\eqref{e_function_sec_1}. An application of Schwarz inequality followed by Lemma~\ref{l_fourth_moment} yields
\begin{align}
T_{\eqref{e_function_sec_1}} &\leq \| f- \mathbb{E}_{\mu^{\sigma}} \left[ f(X)\right] \|_{L^2 (\mu^{\sigma} )} \left( \mathbb{E}_{\mu^{\sigma}} \left[ \left( \sum_{i \in S} \left( X_i -m_i \right) \right)^4 \right] \right)^{\frac{1}{2}} \\
& \lesssim \| f- \mathbb{E}_{\mu^{\sigma}} \left[ f(X)\right] \|_{L^2 (\mu^{\sigma} )} |S| =  \| f- \mathbb{E}_{\mu^{\sigma}} \left[ f(X)\right] \|_{L^2 (\mu^{\sigma} )} |\supp f| \\
& \overset{Poincare}{\lesssim} \| \nabla f \|_{L^2 (\mu^{\sigma} )} |\supp f| \leq \| \nabla f \|_{L^4 (\mu^{\sigma} )} |\supp f|.
\end{align}
Let us turn to the estimation of~\eqref{e_function_sec_2}. As in the proof of Lemma~\ref{l_function_first_estimate}, we denote
\begin{align}
S_d : = \{j \in [N] \ : \ \text{dist}(S, j) =d \}, \qquad d = 1, 2, \cdots.
\end{align}
Let us recall that for each~$d \geq 1$ we have~$|S_d | \leq 2 |\supp f|$.
We write~$T_{\eqref{e_function_sec_2}}$ as
\begin{align}
T_{\eqref{e_function_sec_2}} = 2 \sum_{d \geq 1} \cov_{\mu^{\sigma}} \left( \left( f(X) - \mathbb{E}_{\mu^{\sigma}}\left[ f(X) \right] \right) \sum_{i \in S} \left( X_i -m_i \right)  ,  \sum_{k \in S_d} X_k \right).
\end{align}
It holds by Theorem~\ref{p_decay_of_correlations_gce} that
\begin{align}
&\cov_{\mu^{\sigma}} \left( \left( f(X) - \mathbb{E}_{\mu^{\sigma}}\left[ f(X) \right] \right) \sum_{i \in S} \left( X_i -m_i \right)  ,  \sum_{k \in S_d} X_k \right) \\
& \qquad \lesssim \left\| \nabla \left( \left( f(X) - \mathbb{E}_{\mu^{\sigma}}\left[ f(X) \right] \right) \sum_{i \in S} \left( X_i -m_i \right) \right) \right\|_{L^2 ( \mu^{\sigma}) } |S_d|^{\frac{1}{2}} \exp \left( - C d \right) \\
& \qquad \lesssim \left\| \left( \sum_{i \in S} \left( X_i -m_i \right) \right) \nabla f(X)  \right\|_{L^2 (\mu^{\sigma})} |\supp f|^{\frac{1}{2}} \exp \left( - C d \right) \label{e_func_sec_2_1}\\
& \qquad \quad + \left\| \left( f(X) - \mathbb{E}_{\mu^{\sigma}}\left[ f(X) \right] \right)  \right\|_{L^2 (\mu^{\sigma})}|\supp f| \exp \left( - Cd \right). \label{e_func_sec_2_2}
\end{align}
Then a direct calculation yields
\begin{align}
\left\| \left( \sum_{i \in S} \left( X_i -m_i \right) \right) \nabla f(X)  \right\|_{L^2 (\mu^{\sigma})}^2 & = \mathbb{E}_{\mu^{\sigma}} \left[ \left(\sum_{i \in S} \left( X_i -m_i \right) \right)^2 | \nabla f(X) |^2  \right] \\
& \leq \| \nabla f \|_{L^4 (\mu^{\sigma})}^2  \mathbb{E}_{\mu^{\sigma}} \left[ \left(\sum_{i \in S} \left( X_i -m_i \right) \right)^4 \right]^{\frac{1}{2}} \\
& \overset{Corollary~\ref{l_moment_estimate}}{\lesssim} \| \nabla f \|_{L^4 (\mu^{\sigma})}^2 |S| = \| \nabla f \|_{L^4 (\mu^{\sigma})}^2 |\supp f|. \label{e_func_sec_2_1_est}
\end{align}
Hence plugging the estimate~\eqref{e_func_sec_2_1_est} into~\eqref{e_func_sec_2_1} gives
\begin{align}
T_{\eqref{e_func_sec_2_1}} \lesssim \| \nabla f \|_{L^4 (\mu^{\sigma})}^2 |\supp f| \exp \left( - Cd \right).
\end{align}
It also holds from Poincar\'e inequality that
\begin{align}
T_{\eqref{e_func_sec_2_2}} \lesssim  \| \nabla f \|_{L^2 (\mu^{\sigma})}|\supp f| \exp \left( -Cd\right) \leq \| \nabla f \|_{L^4 (\mu^{\sigma})}|\supp f| \exp \left( -Cd\right).
\end{align}
Thus we conclude
\begin{align}
T_{\eqref{e_function_sec_2}} \lesssim \sum_{ d \geq 1 } \left( T_{\eqref{e_func_sec_2_1}} + T_{\eqref{e_func_sec_2_2}} \right)  \lesssim \sum_{ d \geq 1 } \| \nabla f \|_{L^4 (\mu^{\sigma})}|\supp f| \exp \left( -Cd\right) \lesssim  \| \nabla f \|_{L^4 (\mu^{\sigma})}|\supp f|.
\end{align}
The estimation of~\eqref{e_function_sec_3} follows from similar calculations given in Lemma~\ref{l_third_moment} and estimation of~\eqref{e_function_sec_2}.
\qed
\medskip 

For any~$i \in [N]$ and~$l >0$, let us denote~$B_l (i)$ by
\begin{align}
B_l (i ) : = \{ k \in [N] \ : \ |k-i| \leq l \}. \label{d_def_bli}
\end{align}

\medskip 

The last statement of this section provides a lower bound onto block-block correlations.

\begin{lemma} \label{l_choice_of_l_positive}
There are positive constants~$C$ and~$L$ such that for any~$i \in [N]$,~$l \geq L$, and a set~$T \supseteq B_l (i)$,
\begin{align}
\cov_{\mu^{\sigma}}\left( \sum_{j \in B_l (i) } X_j , \sum_{k \in T} X_k \right) \geq Cl.
\end{align}
\end{lemma}

\medskip 

\noindent \emph{Proof of Lemma~\ref{l_choice_of_l_positive}.} \ By Lemma~\ref{l_variance_estimate}, there is a constant~$C_1 >0$ such that
\begin{align}
\var_{\mu^{\sigma}}\left( \sum_{j \in B_l (i) } X_j  \right) \geq C_1 l.
\end{align}
Next, an application of Theorem~\ref{p_decay_of_correlations_gce} implies that for some~$C_2 >0$,
\begin{align}
\sum_{k \in T \backslash B_l (i)} \left|\cov_{\mu^{\sigma}}\left( \sum_{j \in B_l (i) } X_j ,  X_k \right)  \right| &\leq C\sum_{k \in T \backslash B_l (i)} l^{\frac{1}{2}}\exp \left(- C \text{dist}(k, B_l (i) ) \right) \\
&\leq C l^{\frac{1}{2}} \sum_{d=1}^{\infty} \exp \left( -Cd \right) \leq C_2 l^{\frac{1}{2}}.
\end{align}
We note that the constants~$C_1$ and~$C_2$ are uniform. By choosing~$L$ large enough, there is a uniform constant~$C$ such that for any~$l \geq L$,
\begin{align}
\cov_{\mu^{\sigma}}\left( \sum_{j \in B_l (i) } X_j , \sum_{k \in T} X_k \right) &= \var_{\mu^{\sigma}}\left( \sum_{j \in B_l (i) } X_j  \right) - \sum_{k \in T \backslash B_l (i)} \cov_{\mu^{\sigma}}\left( \sum_{j \in B_l (i) } X_j ,  X_k \right) \\
& \geq C_1 l - C_2 l^{\frac{1}{2}} \geq C l.
\end{align}
\qed

\medskip 

\subsection{Influence of boundary conditions onto observables and correlations.}

In this section,~$f$ is a given intensive function and let us denote~$S = \supp f$. We decompose the sublattice~$[N]$ into two sets~$E_S$ and~$F_S$ as follows:
\begin{align}
E_S : = \{ 1, \cdots, N\} \cap \{ k : dist(k , S ) \leq M \ln N \}, \label{d_def_e1f} \\
F_S : = \{ 1, \cdots, N\} \cap \{ k : dist(k , S ) > M \ln N \}, \label{d_def_e2f}
\end{align}
where~$M$ is a sufficiently large constant which will be chosen later. Recalling the definition~\eqref{d_def_sd} of~$S_d$, the sets~$E_S$ and~$F_S$ can be written as
\begin{align}
E_S = \bigcup_{d=0}^{M \ln N} S_d , \qquad F_S = \bigcup_{d=M \ln N +1}^{N} S_d, 
\end{align}
where we define~$S_0$ to be~$S$. We decompose the gce~$\mu^{\sigma}$ into the conditional measure~$\mu^{\sigma}\left( dx^{E_S} \ \big| \ y^{F_S}  \right)$ and the marginal measure~$\bar{\mu}^{\sigma} \left( dy ^{F_S}  \right)$. That is, for any test function~$\zeta$
\begin{align} \label{e_measure_decomposition}
\int \zeta \mu^{\sigma} = \int \int \zeta \left(x^{E_S} , y^{F_S} \right) \mu^{\sigma}\left( dx^{E_S} \ \big| \ y^{F_S}  \right)\bar{\mu}^{\sigma} \left( dy ^{F_S}  \right).
\end{align}
To reduce our notational burden, we write~$x= x^{E_S},$~$y = y^{F_S}$, and $z= z^{F_S}$ in this section. \\

The next two lemmas estimate the influence of the boundary conditions onto observables and correlations. \\

\begin{lemma} \label{l_conditional_function_difference}
Let~$y = y^{F_S}$ and~$z = z^{F_S}$ be given. For~$N$ large enough, it holds that
\begin{align}
&\left| \mathbb{E}_{\mu^{\sigma}( dx  |  y  )} \left[ f(X) \right] -  \mathbb{E}_{\mu^{\sigma}( dx  |  z  )} \left[ f(X) \right] \right| \\
& \qquad \lesssim \| \nabla f \|_{\infty}  \left( \sum_{\substack{i \in E_S , j \in F_S \\ |i-j| \leq R}} M_{ij}^2(y_j -z_j )^2 \right)^{\frac{1}{2}} \exp \left( -C M \ln N\right).
\end{align}
\end{lemma}

\medskip

\noindent \emph{Proof of Lemma~\ref{l_conditional_function_difference}.} \ By interpolation we have
\begin{align}
&\mathbb{E}_{\mu^{\sigma}(dx|y)} \left[ f(X) \right] - \mathbb{E}_{\mu^{\sigma}(dx|z)} \left[ f(X) \right] \\
& \qquad = \int_0 ^1 \frac{d}{dt} \mathbb{E}_{\mu^{\sigma}(dx|ty + (1-t)z)} \left[ f(X) \right] dt \\
& \qquad = \int_0 ^1 \cov_{\mu^{\sigma}(dx|ty + (1-t)z)} \left( f  ,  \sum_{\substack{i \in E_S , j \in F_S \\ |i-j| \leq R}} M_{ij} X_i (y_j -z_j) \right) dt.  \label{e_interpolation_gce}
\end{align}
With slight abuse of notation we denote~$\mu_{c,t}^{\sigma} = \mu^{\sigma} (dx | ty + (1-t)z)$. We note that the conditional measure~$\mu_{c, t}^{\sigma}$ is again a gce. First, we compute
\begin{align}
\left\| \nabla \left( \sum_{\substack{i \in E_S , j \in F_S \\ |i-j| \leq R}} M_{ij} X_i (y_j -z_j) \right) \right\|_{L^2 (\mu_{c,t}^{\sigma})} ^2 &= \mathbb{E}_{\mu_{c,t}^{\sigma}} \left[ \sum_{i \in E_S} \left( \sum_{\substack{j \in F_S \\ |i-j| \leq R}} M_{ij} (y_j -z_j ) \right)^2 \right]  \\
& =\sum_{i \in E_S} \left( \sum_{\substack{j \in F_S \\ |i-j| \leq R}} M_{ij} (y_j -z_j ) \right)^2 \\
& \leq  \sum_{i \in E_S} (2R) \sum_{\substack{j \in F_S \\ |i-j| \leq R}} M_{ij}^2 (y_j -z_j )^2  \\
& \lesssim  \sum_{\substack{i \in E_S , j \in F_S \\ |i-j| \leq R}} M_{ij}^2(y_j -z_j )^2, \label{e_conditional_auxiliary_computation}
\end{align}
where the first inequality follows from Cauchy's inequality and an observation that for each~$i \in E_S$ there are at most~$2R$ number of~$j$'s with~$j \in F_S$ and~$|i-j| \leq R$. \\

For a pair~$(i,j)$ with~$i \in E_S$,~$j \in F_S$ and~$|i-j| \leq R$, the triangle inequality implies that for large enough~$N$
\begin{align}
\text{dist}(i, S) \geq \text{dist}(j, S ) - |i-j| \geq M \ln N - R \geq \frac{1}{2} M \ln N.
\end{align}
Because~$\mu_{c,t}^{\sigma}$ is also a gce, an application of Theorem~\ref{p_decay_of_correlations_gce} yields the desired estimate for the integrand in~\eqref{e_interpolation_gce}:
\begin{align}
&\left| \cov_{\mu_{c,t}^{\sigma}} \left( f  ,  \sum_{\substack{i \in E_S , j \in F_S \\ |i-j| \leq R}} M_{ij} X_i (y_j -z_j) \right) \right| \\
& \qquad \lesssim  \| \nabla f \|_{L^2 (\mu_{c,t}^{\sigma})}\left\| \nabla \left( \sum_{\substack{i \in E_S , j \in F_S \\ |i-j| \leq R}} M_{ij} X_i (y_j -z_j) \right) \right\|_{L^2 (\mu_{c,t}^{\sigma})} \exp \left( - C \frac{1}{2} M \ln N \right) \\
&\qquad \overset{\eqref{e_conditional_auxiliary_computation}}{\lesssim} \| \nabla f \|_{\infty} \left( \sum_{\substack{i \in E_S , j \in F_S \\ |i-j| \leq R}} M_{ij}^2(y_j -z_j )^2 \right)^{\frac{1}{2}} \exp \left( -C M \ln N\right). \label{e_interpolation_estimate}
\end{align}
Plugging the estimate~\eqref{e_interpolation_estimate} into~\eqref{e_interpolation_gce} finishes the proof of Lemma~\ref{l_conditional_function_difference}.
\qed 

\medskip

\begin{lemma} \label{l_conditional_covariance_difference}
Under the same assumptions as in Lemma~\ref{l_conditional_function_difference}, we have for each~$k \in E_S$ with~$\text{dist}(S, k) \geq \frac{1}{2} M \ln N$, 
\begin{align} \label{e_far_boundary}
&\left|\cov_{\mu^{\sigma}(dx | y ) } \left( f(X)  ,  X_k \right)-\cov_{\mu^{\sigma}(dx | z ) } \left( f(X)  ,  X_k \right) \right| \lesssim \| \nabla f \|_{\infty} \exp \left( -C \text{dist}(S, k) \right) ,
\end{align} 
and for~$k \in E_S$ with~$\text{dist}(S, k)  < \frac{1}{2} M \ln N$,
\begin{align} 
&\left|\cov_{\mu^{\sigma}(dx | y ) } \left( f(X)  ,  X_k \right)-\cov_{\mu^{\sigma}(dx | z ) } \left( f(X)  ,  X_k \right) \right|  \\
& \qquad \lesssim \| \nabla f \|_{\infty} \left( \sum_{\substack{i \in E_S , j \in F_S \\ |i-j| \leq R}} M_{ij}^2(y_j -z_j )^2 \right)^{\frac{1}{2}} \exp \left( -C M \ln N\right). \label{e_near_boundary}
\end{align}
\end{lemma}
\medskip

\begin{remark}
The first estimate~\eqref{e_far_boundary} implies that if the distance between~$S$ and~$k$ are far enough, the covariances are uniformly bounded from above. If~$S$ and~$k$ are close,~\eqref{e_near_boundary} implies that the difference of covariances can be measured in terms of boundary spin values.
\end{remark}
\medskip

\noindent \emph{Proof of Lemma~\ref{l_conditional_covariance_difference}.} \ Because conditional measures~$\mu^{\sigma}(dx|y)$ and~$\mu^{\sigma}(dx|z)$ are again gces, the case when~$\text{dist}(S , k) \geq \frac{1}{2}M \ln N$ directly follows from Theorem~\ref{p_decay_of_correlations_gce}. Let us assume~$\text{dist}(S, k) < \frac{1}{2}M \ln N$. As in the proof of Lemma~\ref{l_conditional_function_difference} we use interpolation to get
\begin{align}
&\cov_{\mu^{\sigma}(dx | y ) } \left( f(X)  ,  X_k \right)-\cov_{\mu^{\sigma}(dx | z ) } \left( f(X)  ,  X_k \right) \\
& \qquad  = \int_0 ^1 \frac{d}{dt} \cov_{\mu_{c,t}^{\sigma}} \left( f(X)  ,  X_k \right) dt,
\end{align}
where~$\mu_{c,t}(dx) = \mu^{\sigma} (dx | ty + (1-t)z)$. A straightforward calculation gives
\begin{align}
&\frac{d}{dt} \cov_{\mu_{c,t}^{\sigma}} \left(f(X) , X_t \right) \\
& = \frac{d}{dt} \mathbb{E}_{\mu_{c,t}^{\sigma}} \left[ \left(f(X) - \mathbb{E}_{\mu_{c,t}^{\sigma}} \left[ f(X) \right] \right)\left( X_k - \mathbb{E}_{\mu_{c,t}^{\sigma}} \left[X_k \right] \right) \right] \\
& = \cov_{\mu_{c,t}^{\sigma}} \left(\left(f(X) - \mathbb{E}_{\mu_{c,t}^{\sigma}} \left[ f(X) \right] \right)\left( X_k - \mathbb{E}_{\mu_{c,t}^{\sigma}} \left[X_k \right] \right)   ,  \sum_{\substack{i \in E_S , j \in F_S \\ |i-j| \leq R}} M_{ij} X_i (y_j -z_j) \right)  \label{e_conditional_covariance_first}  \\
& \quad -  \cov_{\mu_{c,t}^{\sigma}} \left(f(X) - \mathbb{E}_{\mu_{c,t}^{\sigma}} \left[ f(X) \right]   ,  \sum_{\substack{i \in E_S , j \in F_S \\ |i-j| \leq R}} M_{ij} X_i (y_j -z_j) \right) \mathbb{E}_{\mu_{c,t}^{\sigma}} \left[  X_k - \mathbb{E}_{\mu_{c,t}^{\sigma}} \left[X_k \right] \right] \label{e_conditional_covariance_second}\\
& \quad - \cov_{\mu_{c,t}^{\sigma}} \left(   X_k - \mathbb{E}_{\mu_{c,t}^{\sigma}} \left[X_k \right]    ,  \sum_{\substack{i \in E_S , j \in F_S \\ |i-j| \leq R}} M_{ij} X_i (y_j -z_j) \right) \mathbb{E}_{\mu_{c,t}^{\sigma}} \left[f(X) - \mathbb{E}_{\mu_{c,t}^{\sigma}} \left[ f(X) \right] \right]. \label{e_conditional_covariance_third}
\end{align}
Then Theorem~\ref{p_decay_of_correlations_gce}, Corollary~\ref{l_moment_estimate} and Lemma~\ref{l_function_moment_estimate} imply (see also estimations of~\eqref{e_func_sec_2_1},~\eqref{e_func_sec_2_2}, and~\eqref{e_conditional_auxiliary_computation})
\begin{align}
\left| T_{\eqref{e_conditional_covariance_first}} \right| &\lesssim \left( \| \nabla f \|_{\infty}  + \| f(X) - \mathbb{E}_{\mu_{c,t}} \left[ f(X) \right]\|_{L^2(\mu_{c,t}^{\sigma})} \right) \\
& \qquad \times \left( \sum_{\substack{i \in E_S , j \in F_S \\ |i-j| \leq R}} M_{ij}^2 (y_j -z_j )^2 \right)^{\frac{1}{2}} \exp \left( -C M \ln N\right) \\
& \lesssim \| \nabla f \|_{\infty} \left( \sum_{\substack{i \in E_S , j \in F_S \\ |i-j| \leq R}} M_{ij}^2(y_j -z_j )^2 \right)^{\frac{1}{2}} \exp \left( -C M \ln N\right).
\end{align}
Similar calculation also yields
\begin{align}
\left| T_{\eqref{e_conditional_covariance_second}} \right|, \left| T_{\eqref{e_conditional_covariance_third}} \right| \lesssim \| \nabla f \|_{\infty} \left( \sum_{\substack{i \in E_S , j \in F_S \\ |i-j| \leq R}} M_{ij}^2(y_j -z_j )^2 \right)^{\frac{1}{2}} \exp \left( -C M \ln N\right).
\end{align}
Hence we get
\begin{align}
&\left| \cov_{\mu^{\sigma}(dx | y ) } \left( f(X) \ , \ X_k \right)-\cov_{\mu^{\sigma}(dx | z ) } \left( f(X) \ , \ X_k \right) \right| \\
& \qquad \lesssim \int_{0}^{1}  \| \nabla f \|_{\infty}  \left( \sum_{\substack{i \in E_S , j \in F_S \\ |i-j| \leq R}} M_{ij}^2(y_j -z_j )^2 \right)^{\frac{1}{2}} \exp \left( -C M \ln N\right) dt \\ 
& \qquad =  \| \nabla f \|_{\infty}  \left( \sum_{\substack{i \in E_S , j \in F_S \\ |i-j| \leq R}} M_{ij}^2(y_j -z_j )^2 \right)^{\frac{1}{2}} \exp \left( -C M \ln N\right).
\end{align}
\qed

\medskip

\begin{corollary} \label{l_conditional_covariance_corollary}
Under the same assumptions as in Lemma~\ref{l_conditional_function_difference}, we have for each~$k \in E_S$ with~$\text{dist}(S, k) \geq \frac{1}{2} M \ln N$, 
\begin{align}
&\left|\cov_{\mu^{\sigma}(dx | y ) } \left( f(X) \ , \ X_k \right)-\cov_{\mu^{\sigma} } \left( f(X) \ , \ X_k \right) \right|  \lesssim \| \nabla f \|_{\infty} \exp \left( -C \text{dist}(S, k) \right) ,
\end{align}
and for~$k \in E_S$ with~$\text{dist}(S, k)  < \frac{1}{2} M \ln N$,
\begin{align}
&\left|\cov_{\mu^{\sigma}(dx | y ) } \left( f(X) \ , \ X_k \right)-\cov_{\mu^{\sigma}} \left( f(X) \ , \ X_k \right) \right|  \\
& \qquad \lesssim \| \nabla f \|_{\infty} \exp \left( -C M \ln N\right) \int \left( \sum_{\substack{i \in E_S , j \in F_S \\ |i-j| \leq R}} M_{ij}^2(y_j -z_j )^2 \right)^{\frac{1}{2}} \bar{\mu}^{\sigma} (dz) \\
& \qquad \quad + \| \nabla f \|_{\infty} |\supp f |^{\frac{1}{2}} \exp \left( - C M \ln N\right).
\end{align}
\end{corollary}

\medskip

\noindent \emph{Proof of Corollary~\ref{l_conditional_covariance_corollary}.} \ The first case follows from Theorem~\ref{p_decay_of_correlations_gce} and triangle inequality. To prove the case when~$k \in E_S$,~$\text{dist}(S, k) < \frac{1}{2} M \ln N$, we use the law of total covariance and write
\begin{align}
&\left|\cov_{\mu^{\sigma}(dx | y ) } \left( f(X)  ,  X_k \right)-\cov_{\mu^{\sigma}} \left( f(X)  ,  X_k \right) \right|\\
& \qquad \leq \left|\cov_{\mu^{\sigma}(dx | y ) } \left( f(X)  ,  X_k \right) - \int \cov_{\mu^{\sigma}(dx | z ) } \left( f(X)  ,  X_k \right) \bar{\mu}^{\sigma} (dz) \right|  \label{e_cond_cov_cor_eq1}\\
& \qquad \quad + \left| \cov_{\mu^{\sigma}} \left( \mathbb{E}_{\mu^{\sigma}(dx|y)} \left[f(X) \right] \ ,\  \mathbb{E}_{\mu^{\sigma}(dx|y)} \left[ X_k \right] \right) \right| . \label{e_cond_cov_cor_eq2}
\end{align}
Then Lemma~\ref{l_conditional_covariance_difference} implies
\begin{align}
T_{\eqref{e_cond_cov_cor_eq1}} &\leq \int \left| \cov_{\mu^{\sigma}(dx | y ) } \left( f(X)  ,  X_k \right) -  \cov_{\mu^{\sigma}(dx | z ) } \left( f(X)  ,  X_k \right)  \right| \bar{\mu}^{\sigma} (dz) \\
& \leq \| \nabla f \|_{\infty} \exp \left( -C M \ln N\right) \int \left( \sum_{\substack{i \in E_S , j \in F_S \\ |i-j| \leq R}} M_{ij}^2(y_j -z_j )^2 \right)^{\frac{1}{2}} \bar{\mu}^{\sigma} (dz).
\end{align}
A similar calculations using Lemma~\ref{l_conditional_function_difference} gives
\begin{align}
T_{\eqref{e_cond_cov_cor_eq2}} &\leq \int \int \left| \mathbb{E}_{\mu^{\sigma}(dx|y)} \left[f(X) \right] - \mathbb{E}_{\mu^{\sigma}(dx|z)} \left[f(X) \right] \right| \\
& \qquad \qquad \times \left| \mathbb{E}_{\mu^{\sigma}(dx|y)} \left[X_k \right] - \mathbb{E}_{\mu^{\sigma}(dx|z)} \left[X_k \right] \right| \bar{\mu}^{\sigma} (dy)\bar{\mu}^{\sigma} (dz) \\
& \lesssim \| \nabla f \|_{\infty}  \exp \left( - 2C M \ln N\right) \\
&\quad \times \int \int \left( \sum_{\substack{i \in E_S , j \in F_S \\ |i-j| \leq R}} M_{ij}^2(y_j -z_j )^2 \right)^{\frac{1}{2}} \left( \sum_{\substack{i \in E_k , j \in F_k \\ |i-j| \leq R}} M_{ij}^2(y_j -z_j )^2 \right)^{\frac{1}{2}}\bar{\mu}^{\sigma} (dy)\bar{\mu}^{\sigma} (dz),
\end{align}
where~$E_k$ and~$F_k$ are defined by
\begin{align}
E_k : = \{ 1, \cdots, N\} \cap \{ l : dist(l , k ) \leq M \ln N \}, \\
F_k : = \{ 1, \cdots, N\} \cap \{ l : dist(l , k ) > M \ln N \}.
\end{align}
Now we apply Schwarz inequality followed by Corollary~\ref{l_moment_estimate} and get, as desired,
\begin{align}
&\int \int \left( \sum_{\substack{i \in E_S , j \in F_S \\ |i-j| \leq R}} M_{ij}^2(y_j -z_j )^2 \right)^{\frac{1}{2}} \left( \sum_{\substack{i \in E_k , j \in F_k \\ |i-j| \leq R}} M_{ij}^2(y_j -z_j )^2 \right)^{\frac{1}{2}}\bar{\mu}^{\sigma} (dy)\bar{\mu}^{\sigma} (dz) \\
&\qquad \leq \left(\int \int  \sum_{\substack{i \in E_S , j \in F_S \\ |i-j| \leq R}} M_{ij}^2(y_j -z_j )^2  \bar{\mu}^{\sigma} (dy)\bar{\mu}^{\sigma} (dz) \right)^{\frac{1}{2}} \\
&\qquad \quad \times \left(\int \int  \sum_{\substack{i \in E_k , j \in F_k \\ |i-j| \leq R}} M_{ij}^2(y_j -z_j )^2  \bar{\mu}^{\sigma} (dy)\bar{\mu}^{\sigma} (dz) \right)^{\frac{1}{2}} \\
&\qquad  \lesssim \left(\sum_{\substack{i \in E_S , j \in F_S \\ |i-j| \leq R}} M_{ij}^2 \var\left( X_j \right) \right)^{\frac{1}{2}}  \left(\sum_{\substack{i \in E_k , j \in F_k \\ |i-j| \leq R}} M_{ij}^2 \var\left( X_j \right) \right)^{\frac{1}{2}} \\
& \qquad \lesssim \left(2R^2 |\supp f| \right)^{\frac{1}{2}} \left(2R^2 \right)^{\frac{1}{2}} \sim |\supp f|^{\frac{1}{2}}.
\end{align}
\qed 
\medskip

\section{Proof of Theorem~\ref{p_equivalence_observables} and Theorem~\ref{p_equivalence_decay_of_correlations}} \label{s_proof_equiv_obs}

As it is common when deducing equivalence of ensembles, we express the difference of observables and correlations between gce and ce using the inverse Fourier transform. \\

\begin{lemma} \label{l_inverse_fourier}
For any function $\zeta, \eta : \mathbb{R}^N \to \mathbb{R}$,
\begin{align} 
&\mathbb{E}_{\mu_m} \left[ \zeta (X) \right] - \mathbb{E}_{\mu^{\sigma}} \left[ \zeta (X) \right]\\
& \qquad = \frac{\int_{\mathbb{R}}\mathbb{E}_{\mu^{\sigma}} \left[ \left( \zeta (X) - \mathbb{E}_{\mu^{\sigma}} \left[ \zeta (X) \right]\right) \exp \left( i \frac{1}{\sqrt{N}} \sum_{i=1}^{N} \left(X_i -m_i \right) \xi\right) \right]d\xi}{\int_{\mathbb{R}}\mathbb{E}_{\mu^{\sigma}} \left[  \exp \left( i \frac{1}{\sqrt{N}} \sum_{i=1}^{N} \left(X_i -m_i \right) \xi\right) \right]d\xi}
\end{align}
and 
\begin{align}
&\cov_{\mu_m} \left( \zeta(X)  ,  \eta (X) \right) - \cov_{\mu^{\sigma}} \left( \zeta (X)  ,  \eta(X) \right) \\
& = \frac{\int_{\mathbb{R}}\mathbb{E}_{\mu^{\sigma}} \left[ \left( \zeta (X) - \mathbb{E}_{\mu^{\sigma}} \left[ \zeta (X) \right]\right)\left( \eta (X) - \mathbb{E}_{\mu^{\sigma}} \left[ \eta (X) \right]\right) \exp \left( i \frac{1}{\sqrt{N}} \sum_{i=1}^{N} \left(X_i -m_i \right) \xi\right) \right]d\xi}{\int_{\mathbb{R}}\mathbb{E}_{\mu^{\sigma}} \left[  \exp \left( i \frac{1}{\sqrt{N}} \sum_{i=1}^{N} \left(X_i -m_i \right) \xi\right) \right]d\xi} \\
& \quad - \left(\mathbb{E}_{\mu_m} \left[ \zeta (X) \right] - \mathbb{E}_{\mu^{\sigma}} \left[ \zeta (X) \right]\right) \left(\mathbb{E}_{\mu_m} \left[ \eta (X) \right] - \mathbb{E}_{\mu^{\sigma}} \left[ \eta (X) \right]\right).
\end{align}
\end{lemma}

\medskip 

The proof of Lemma~\ref{l_inverse_fourier} is outlined in the Appendix. We provide the proof of Theorem~\ref{p_equivalence_observables} and Theorem~\ref{p_equivalence_decay_of_correlations} in Section~\ref{s_proof_1st_thm} and Section~\ref{s_proof_2nd_thm}, respectively. \\

\subsection{Proof of Theorem~\ref{p_equivalence_observables}} \label{s_proof_1st_thm}

The proof of Theorem~\ref{p_equivalence_observables} is quite technical. In a naive approach of~\cite{kwme18a}, a similar result was deduced with sub-optimal scaling on~$N$. The main technical difficulty comes from the estimation of the first order term in a Taylor expansion. Let us outline how we overcome this obstacle. \\

We decompose the sublattice~$[N]$ into large blocks (cf. Lemma~\ref{l_choice_of_l_positive}). Then for each intensive observable~$f$, we carefully choose a linear approximation~$h_f$ in terms of block spins (cf.~\eqref{e_def_hf}). The key observation is that the difference~$h = f- h_f$ satisfies the equivalence of observables of the right order because, by choosing~$h_f$ wisely, the problematic first order term in the Taylor expansion becomes small (see proof of Proposition~\ref{p_1st_der_main_proposition} for more details). We then show that each summand in~$h_f$ also satisfies the equivalence of ensembles with the right scaling (cf. Lemma~\ref{l_spin_auxiliary_estimates}). This step takes advantage of the fact that~$h_f$ is linear and that the gce and ce have the same mean. Hence, together with elementary estimations, the linear function~$h_f$ can be decomposed into summands sharing a similar structure as~$h= f-h_f$ (cf.~\eqref{e_def_hij}); and therefore also satisfy the equivalence of observables of the right order. We refer to Lemma~\ref{l_spin_auxiliary_estimates} for more details. \\

Now let us turn to the detailed arguments. Let us begin with introducing auxiliary notations and definitions that are needed for the proof of Theorem~\ref{p_equivalence_observables}. Recalling the definition~\eqref{d_def_bli} of~$B_l (i)$, we decompose~$[N]$ as
\begin{align} \label{e_decomposition_of_N}
[N] = \bigcup_{j=1}^{M} B_{l_j}(w_j),
\end{align}
where~$w_j \in [N]$,~$l_j \geq L$ as in Lemma~\ref{l_choice_of_l_positive}, and the union is disjoint. For notational simplicity, we denote for each~$j \in [M]$,
\begin{align} \label{e_decomp_ball_notation}
B_j = B_{l_j} (w_j).
\end{align}
Then the decomposition~\eqref{e_decomposition_of_N} is rewritten as
\begin{align} \label{e_simplified_decomposition}
    [N] = \bigcup_{j=1}^M B_j .
\end{align}

Let us define a map~$\varphi : [N] \to [M]$ that matches each site~$i \in [N]$ with the block that contains it. More precisely, for each~$i \in [N]$, there exists a unique~$j(i) \in [M]$ such that~$i \in B_{j(i)}$. Let us write
\begin{align}
     B_{\varphi(i)} = B_{j(i)} = B_{l _{j(i)}} (w_{j(i)}) \qquad \text{for each } i \in [N].
\end{align}

\medskip
The first step towards to the proof of Theorem~\ref{p_equivalence_observables} is considering a special form of functions. Let us
recall the definitions~\eqref{d_def_e1f} and~\eqref{d_def_e2f} of~$E_S$ and~$F_S$, respectively. Let us fix an intensive function~$f$. We denote~$S = \supp f$ and define
\begin{align} \label{e_def_cf}
c_f : = \frac{\cov_{\mu^{\sigma}} \left( f(X)  ,  \sum_{j \in E_S } X_j  \right)}{\cov_{\mu^{\sigma}}\left( \sum_{i \in S} \sum_{k \in B_{\varphi(i)}} X_k ,  \sum_{j \in E_S} X_j  \right)}.
\end{align}
\medskip 
\begin{remark} \label{r_denom_positive}
By choosing~$l$ large enough, the denominator of~$c_f$ is bounded from below and hence~$c_f$ is well defined. More precisely, by Lemma~\ref{l_choice_of_l_positive} we have
\begin{align}
 \cov_{\mu^{\sigma}}\left( \sum_{i \in S}\sum_{k \in B_{\varphi(i)}} X_k  , \sum_{j \in E_S} X_j  \right) \geq CL |S| \gtrsim |\supp f| > 0.
\end{align}
Moreover, combined with Lemma~\ref{l_function_first_estimate} we have the following estimate:
\begin{align} \label{e_alpha_g_estimate} 
\left| c_f \right| \lesssim \frac{\| \nabla f \|_{\infty} }{|\supp f |^{\frac{1}{2}}}.
\end{align}
\end{remark}

\medskip 

We define a linear approximation~$h_f$ of~$f$ as
\begin{align}\label{e_def_hf}
h_f(x) = c_f \sum_{i \in S}\sum_{k \in B_{\varphi(i)}} X_k 
\end{align}
and write the difference~$h$ as
\begin{align} \label{e_assumption_f_structure}
h(x) = f(x) - h_f(x) = f(x) - c_f \sum_{i \in S}\sum_{k \in B_{\varphi(i)}} X_k.
\end{align}

The following proposition contains core estimate needed for the proof of Theorem~\ref{p_equivalence_observables}. \\

\begin{proposition} \label{p_1st_der_main_proposition}
There exist uniform constants~$N_0 \in \mathbb{N}$ and~$C> 0$ independent of the external field~$s$ and the mean spin~$m$ such that for all~$N \geq N_0$,
\begin{align} \label{e_special_form_proposition}
\left|\int_{\mathbb{R}}\mathbb{E}_{\mu^{\sigma}} \left[ \left( h (X) - \mathbb{E}_{\mu^{\sigma}} \left[ h (X) \right]\right) \exp \left( i \frac{1}{\sqrt{N}} \sum_{i=1}^{N} \left(X_i -m_i \right) \xi\right) \right]d\xi \right| \leq C \frac{|\supp f|}{N} \| \nabla f \|_{\infty}.
\end{align}
\end{proposition}

\medskip

\begin{remark} Proposition~\ref{p_1st_der_main_proposition} was motivated by~\cite[Lemma 4.2]{MR00}. The main difference is that~\cite{MR00} considers l-support while the definition of support of a function~$f$ in this paper is the minimal subset of~$\mathbb{Z}$ with~$f(x) = f(x^{\text{supp} f})$. In~\cite{MR00}, the assumption that~$l$ is large enough was used to guarantee the positiveness of~$c_f$ (cf.~\cite[Section 4]{MR00} and Remark~\ref{r_denom_positive}). To address this difference, we artificially introduce the block decomposition~\eqref{e_decomposition_of_N} of~$[N]$ and additionally include a block summation in the definition of~$c_f$ (and consequently~$h(x)$).
\end{remark}

\begin{remark}
One should compare Proposition~\ref{p_1st_der_main_proposition} with~\cite[Proposition 2]{kwme18a}. There, a similar estimate was deduced but scaling on~$N$ is sub-optimal. To improve our estimate, we first consider a special form of functions. One benefit of of considering such functions is that when we apply Taylor expansions to the left hand side of~\eqref{e_special_form_proposition}, the first order term is estimated relatively easily. For more details, we refer to Section~\ref{s_proof_1st_der_main_proposition}.
\end{remark}

\medskip 
We present the proof of Proposition~\ref{p_1st_der_main_proposition} in Section~\ref{s_proof_1st_der_main_proposition}. The following is a direct consequence of Lemma~\ref{l_inverse_fourier}, Proposition~\ref{p_main_computation} and Proposition~\ref{p_1st_der_main_proposition}. \\

\begin{corollary} \label{p_corollary_prop1}
There exist uniform constants~$N_0 \in \mathbb{N}$ and~$C>0$ independent of the external field~$s$ and the mean spin~$m$ such that for all~$N \geq N_0$
\begin{align}
\left| \mathbb{E}_{\mu_m} \left[ h(X) \right] - \mathbb{E}_{\mu^{\sigma}} \left[ h(X) \right]  \right| \leq C \frac{ |\supp f |}{N} \| \nabla f \|_{\infty}.
\end{align}
\end{corollary}
\medskip

We then prove Theorem~\ref{p_equivalence_observables} for~$f = \sum_{k \in B_{\varphi(i)}} X_k $.  \\ 
\begin{lemma} \label{l_spin_auxiliary_estimates}
For each~$i \in S$, it holds that
\begin{align}
\left| \mathbb{E}_{\mu_m} \left[ \sum_{k \in B_{\varphi(i)}} X_k \right] - \mathbb{E}_{\mu^{\sigma}} \left[ \sum_{k \in B_{\varphi(i)}} X_k \right]\right| \lesssim \frac{1}{N}.
\end{align}
\end{lemma}
\medskip

\noindent \emph{Proof of Lemma~\ref{l_spin_auxiliary_estimates}.} \ Let us fix~$i \in S$ and recall the definition~\eqref{e_decomp_ball_notation} of~$B_j$ and the decomposition~\eqref{e_simplified_decomposition}. For each~$j \in [M]$, we set~$S_{ij} : = B_{\varphi(i)} \cup B_j$. Recalling the definition~\eqref{d_def_e1f} of~$E_S$, we analogously denote~$E_{ij}$ and~$F_{ij}$ by
\begin{align}
    E_{ij} : &= [N] \cap \{ k : \text{dist} (k, S_{ij}) \leq M \ln N \}, \\
    F_{ij} : &= [N] \cap \{ k : \text{dist} (k, S_{ij}) > M \ln N \}.
\end{align}

Similar to the way we defined~$h$ from~$f$, let us construct an auxiliary function~$h_{ij}$ from~$\sum_{k \in B_j} X_j$ as follows:
\begin{align} \label{e_def_hij}
    h_{ij} (X) : = \sum_{k \in B_j} X_k - c_{ij} \sum_{k \in B_{\varphi(i)}} X_k,
\end{align}
where
\begin{align}
    c_{ij} = \frac{ \cov_{\mu^{\sigma}} \left( \sum_{k \in B_j} X_k , \sum_{k \in E_{ij }} X_k \right) } { \cov_{\mu^{\sigma}} \left( \sum_{k \in B_{\varphi(i)}} X_k , \sum_{k \in E_{ij}} X_k \right)     }. 
\end{align}
A similar argument using Lemma~\ref{l_choice_of_l_positive} implies that the denominator of~$c_{ij}$ is positive, hence~$h_{ij}$ well defined (cf. Remark~\ref{r_denom_positive}). Moreover, there is a positive constant~$C$ such that
\begin{align} \label{e_positivity_cij}
    \frac{1}{C} \leq c_{ij} \leq C.
\end{align}
\medskip 

A detailed analysis of the proofs show that the arguments for Proposition~\ref{p_1st_der_main_proposition} still apply to~$h_{ij}$ (and hence Corollary~\ref{p_corollary_prop1}) which implies (see Remark~\ref{r_hij} in Section~\ref{s_proof_1st_der_main_proposition} for more details)
\begin{align}
\left| \mathbb{E}_{\mu_m} \left[ h_{ij}(X) \right] - \mathbb{E}_{\mu^{\sigma}} \left[ h_{ij}(X) \right]  \right| \leq C \frac{ 1}{N}.
\end{align}

We additionally observe that (cf.~\eqref{e_m=sum_mi})
\begin{align} \label{e_mean_conservation}
\sum_{k=1}^{N} \mathbb{E}_{\mu_m} \left[ X_k \right] = Nm = \sum_{k=1}^{N} \mathbb{E}_{\mu^{\sigma}} \left[ X_k \right].
\end{align}
Thus an application of triangle inequality yields
\begin{align}
\frac{M}{N} &\gtrsim \left| \sum_{j=1}^M \left( \mathbb{E}_{\mu_m} \left[h_{ij} (X) \right] - \mathbb{E}_{\mu^{\sigma}} \left[ h_{ij} (X) \right] \right) \right| \\
& = \left| \left( \sum_{j=1}^{M} \sum_{k \in B_j} \mathbb{E}_{\mu_m} \left[ X_k \right] - \sum_{j=1}^{M} c_{ij} \sum_{k \in B_{\varphi(i)}} \mathbb{E}_{\mu_m} \left[ X_k \right]  \right) \right. \\
& \qquad \qquad \left.  -  \left( \sum_{j=1}^{M} \sum_{k \in B_j} \mathbb{E}_{\mu^{\sigma}} \left[ X_k \right] - \sum_{j=1}^{M} c_{ij} \sum_{k \in B_{\varphi(i)}} \mathbb{E}_{\mu^{\sigma}} \left[ X_k \right]  \right) \right| \\
& \overset{\eqref{e_simplified_decomposition}, \eqref{e_mean_conservation}}{=} \left|\left( \mathbb{E}_{\mu_m} \left[ \sum_{k \in B_{\varphi(i)}} X_k \right] - \mathbb{E}_{\mu^{\sigma}} \left[ \sum_{k \in B_{\varphi(i)}} X_k \right] \right) \sum_{j =1}^{M} c_{ij}  \right| \\
& \overset{\eqref{e_positivity_cij}}{=}  \left|  \mathbb{E}_{\mu_m} \left[ \sum_{k \in B_{\varphi(i)}} X_k \right] - \mathbb{E}_{\mu^{\sigma}} \left[ \sum_{k \in B_{\varphi(i)}} X_k \right] \right|  \sum_{j=1}^{M} c_{ij} .
\end{align}
Now we conclude from~\eqref{e_positivity_cij} that, as desired,
\begin{align}
\left| \mathbb{E}_{\mu_m} \left[ \sum_{k \in B_{\varphi(i)}} X_k \right] - \mathbb{E}_{\mu^{\sigma}} \left[ \sum_{k \in B_{\varphi(i)}} X_k \right] \right| \lesssim \frac{1}{N}.
\end{align}
\qed

We are now ready to provide the proof of our first main result, Theorem~\ref{p_equivalence_observables}. \\

\noindent \emph{Proof of Theorem~\ref{p_equivalence_observables}.} \ Let us recall the definition~\eqref{e_assumption_f_structure} of~$h$. A combination of Corollary~\ref{p_corollary_prop1}, Lemma~\ref{l_spin_auxiliary_estimates}, and~\eqref{e_alpha_g_estimate} gives
\begin{align}
&\left|\mathbb{E}_{\mu_m} \left[ f(X) \right] - \mathbb{E}_{\mu^{\sigma}} \left[ f(X) \right] \right| \\
&\qquad  \leq \left|\mathbb{E}_{\mu_m} \left[ h(X) \right] - \mathbb{E}_{\mu^{\sigma}} \left[ h(X) \right] \right| + |c_f| \sum_{i \in S}\left|\mathbb{E}_{\mu_m} \left[ \sum_{k \in B_{\varphi(i)}} X_k \right] - \mathbb{E}_{\mu^{\sigma}} \left[ \sum_{k \in B_{\varphi(i)}} X_k \right] \right| \\
&\qquad  \lesssim \frac{|\supp f |}{N} \| \nabla f \|_{\infty} + \frac{\| \nabla f \|_{\infty} }{|\supp f |^{\frac{1}{2}}} |\supp f | \frac{1}{N} \\
& \qquad \lesssim \frac{|\supp f |}{N} \| \nabla f \|_{\infty}.
\end{align}
\qed 

\medskip

\subsection{Proof of Theorem~\ref{p_equivalence_decay_of_correlations}} \label{s_proof_2nd_thm}

The next proposition provides a core estimate that is needed in the proof of Theorem~\ref{p_equivalence_decay_of_correlations}. \\

\begin{proposition} \label{p_2nd_der_main_proposition}
For any intensive functions~$f, g : \mathbb{R}^N \to \mathbb{R}$, there exist constants~$N_0 \in \mathbb{N}$ and $C>0 $ independent of the external field~$s$ and the mean spin~$m$ such that for all~$N\geq N_0$,
\begin{align}
&\left|\int  \mathbb{E}_{\mu^{\sigma}} \left[ \left(f(X) - \mathbb{E}_{\mu^{\sigma}} \left[ f(X)\right] \right)\left(g(X) - \mathbb{E}_{\mu^{\sigma}} \left[ g(X)\right] \right) \exp\left( i \frac{1}{\sqrt{N}} \sum_{k=1}^N \left(X_k -m_k \right) \xi \right) \right] d\xi\right|\\
&\qquad \leq C \| \nabla f \|_{\infty} \| \nabla g \|_{\infty} \left( \frac{|\supp f | + |\supp g|}{N} + \exp \left( -\text{dist} \left( \supp f \ , \ \supp g \right) \right)  \right). \label{e_2nd_derivative_lemma} 
\end{align}
\end{proposition}

\medskip
\begin{remark}
Proposition~\ref{p_2nd_der_main_proposition} is an extension of~\cite[Proposition 3]{kwme18a}. In~\cite{kwme18a}, the authors estimated the left hand side of~\eqref{e_2nd_derivative_lemma} via second order Taylor expansion. In this article, we use third order Taylor expansion combined with fourth moment bounds (Lemma~\ref{l_fourth_moment}) to improve the estimate. For more details, we refer to Section~\ref{s_proof_2nd_der_main_proposition}.
\end{remark}

\medskip

We present the proof of Proposition~\ref{p_2nd_der_main_proposition} in Section~\ref{s_proof_2nd_der_main_proposition}. Let us now provide the proof of Theorem~\ref{p_equivalence_decay_of_correlations}. \\

\noindent \emph{Proof of Theorem~\ref{p_equivalence_decay_of_correlations}.} \ A combination of Lemma~\ref{l_inverse_fourier}, Theorem~\ref{p_equivalence_observables}, Proposition~\ref{p_main_computation} and Proposition~\ref{p_2nd_der_main_proposition} implies that
\begin{align}
&\left| \cov_{\mu^{\sigma}} \left( f(X)  ,  g(X) \right) - \cov_{\mu_m} \left( f (X)  ,  g(X) \right) \right| \\
& \leq \left| \frac{\int_{\mathbb{R}}\mathbb{E}_{\mu^{\sigma}} \left[ \left( f (X) - \mathbb{E}_{\mu^{\sigma}} \left[ f (X) \right]\right)\left( g (X) - \mathbb{E}_{\mu^{\sigma}} \left[ g (X) \right]\right) \exp \left( i \frac{1}{\sqrt{N}} \sum_{i=1}^{N} \left(X_i -m_i \right) \xi\right) \right]d\xi}{\int_{\mathbb{R}}\mathbb{E}_{\mu^{\sigma}} \left[  \exp \left( i \frac{1}{\sqrt{N}} \sum_{i=1}^{N} \left(X_i -m_i \right) \xi\right) \right]d\xi}  \right| \\
& \quad +  \left| \mathbb{E}_{\mu^{\sigma}} \left[ f (X) \right] - \mathbb{E}_{\mu^{\sigma}} \left[ f (X) \right]\right| \left|\mathbb{E}_{\mu^{\sigma}} \left[ g (X) \right] - \mathbb{E}_{\mu^{\sigma}} \left[ g (X) \right]\right| \\
& \lesssim \| \nabla f \|_{\infty} \| \nabla g \|_{\infty} \left( \frac{|\supp f | + |\supp g| }{N} + \exp \left( -\text{dist} \left( \supp f \ , \ \supp g \right) \right)  \right)  \\
& \quad + \frac{|\supp f| |\supp g|}{N^2} \| \nabla f \|_{\infty} \| \nabla g \|_{\infty} \\
& \lesssim \| \nabla f \|_{\infty} \| \nabla g \|_{\infty} \left( \frac{|\supp f | + |\supp g| }{N} + \exp \left( -\text{dist} \left( \supp f \ , \ \supp g \right) \right)  \right). 
\end{align}
\qed

\medskip

\section{Proof of Proposition~\ref{p_1st_der_main_proposition}} \label{s_proof_1st_der_main_proposition}

The main argument for the proof of Proposition~\ref{p_1st_der_main_proposition} follows a well known method for deducing local CLT. Like in the proof of~\cite[Proposition 1]{KwMe18}, the integral is divided into inner and outer parts that are estimated separately. More precisely, let us fix~$\delta>0$ small enough and decompose the integral as follows:
\begin{align}
&\int_{\mathbb{R}} \mathbb{E}_{\mu ^{\sigma} } \left[ \left(h(X)-\mathbb{E}_{\mu^{\sigma}} \left[h(X) \right] \right) \exp\left( i\frac{1}{\sqrt{N}} \sum_{k=1}^N \left(X_k -m_k \right) \xi \right) \right] d\xi \\
& =\int_{\{ \left| \left(1/ \sqrt{N}\right)\xi\right|\leq\delta\}} \mathbb{E}_{\mu^{\sigma}  } \left[ \left(h(X)-\mathbb{E}_{\mu^{\sigma}} \left[h(X) \right] \right) \exp\left( i\frac{1}{\sqrt{N}} \sum_{k=1}^N \left(X_k -m_k \right) \xi \right) \right] d\xi \label{e_exp_decay_inner}  \\
& \quad + \int_{\{ \left| \left(1/ \sqrt{N}\right)\xi\right| > \delta\}} \mathbb{E}_{\mu^{\sigma} } \left[\left(h(X)-\mathbb{E}_{\mu^{\sigma}} \left[h(X) \right] \right) \exp\left( i\frac{1}{\sqrt{N}} \sum_{k=1}^N \left(X_k -m_k \right) \xi \right) \right] d\xi. \label{e_exp_decay_outer}
\end{align}

The estimation of the outer integral~$T_{\eqref{e_exp_decay_outer}}$ is the easy part. In~\cite[Lemma 9]{kwme18a} it was shown that
\begin{align}
    \left| T_{\eqref{e_exp_decay_outer}} \right| \lesssim \|\nabla f \|_{L^2 (\mu^{\sigma})} \frac{ |\supp f |^{\frac{1}{2}}}{N^{\frac{1}{2} - \varepsilon}}.
\end{align}
However, this estimate was stated sub-optimal. As it is usual when deducing local CLTs, the outer integral~$T_{\eqref{e_exp_decay_outer}}$ actually decays exponentially in the system size (see the proof of~\cite[Lemma 9]{kwme18a}). Hence, it holds
\begin{align} \label{e_observable_outer_estimate}
\left| T_{\eqref{e_exp_decay_outer}} \right| \lesssim \|\nabla f \|_{\infty} \frac{ |\supp f |}{N}.
\end{align}

The subtle part of the argument is the estimation of the inner integral~\eqref{e_exp_decay_inner}. \\

\begin{lemma} \label{l_1st_der_inner} It holds that
\begin{align}
&\left|T_{\eqref{e_exp_decay_inner}} \right| \lesssim \| \nabla f \|_{\infty} \frac{ |\supp f |}{N}.
\end{align}
\end{lemma}
\medskip 

\begin{remark}
The estimate of Lemma~\ref{l_1st_der_inner} improves the estimate of~\cite{kwme18a} by a factor of~$\frac{|\supp f|^{\frac{1}{2}}}{N^{\frac{1}{2} + \varepsilon}}$. Essentially, this is due to the special form of~$h$ (see~\eqref{e_assumption_f_structure}) which introduces a quasi-cancellation in a first order Taylor term. For details see proof of Lemma~\ref{l_main_equiv_ob_second_aux}. 
\end{remark}
\medskip
Proposition~\ref{p_1st_der_main_proposition} is a direct consequence of~\eqref{e_observable_outer_estimate} and Lemma~\ref{l_1st_der_inner}. \\

\noindent \emph{Proof of Proposition~\ref{p_1st_der_main_proposition}.} \ A combination of~\eqref{e_observable_outer_estimate} and Lemma~\ref{l_1st_der_inner} proves the Proposition~\ref{p_1st_der_main_proposition}.
\qed

\medskip

 Let us see how we estimate the inner integral~\eqref{e_exp_decay_inner}. We begin with introducing auxiliary definitions and notations for proof of Lemma~\ref{l_1st_der_inner}. We set~$S = \supp f$ and let us recall the definition~\eqref{d_def_e1f} and~\eqref{d_def_e2f} of the sets~$E_S$ and~$F_S$ and the decomposition~\eqref{e_measure_decomposition} of the gce~$\mu^{\sigma}$. To reduce the notational burden we write
\begin{align}
\mu_c ^{\sigma} (dx| y ) = \mu^{\sigma}( dx^{E_S} | y^{F_S}) \qquad \text{and} \qquad  \bar{\mu}^{\sigma} (dy ) = \bar{\mu}^{\sigma} ( dy^{F_S}).
\end{align}
We observe that by the law of total covariance, the integrand in~\eqref{e_exp_decay_inner} can be written as
\begin{align}
&\mathbb{E}_{\mu^{\sigma}} \left[ \left( h(X) - \mathbb{E}_{\mu^{\sigma}} \left[ h(X) \right] \right) \exp \left( i \frac{1}{\sqrt{N}} \sum_{k=1}^{N} \left( X_k -m_k \right) \xi \right) \right] \\
&\qquad = \cov_{\mu^{\sigma}} \left( h(X)  ,    \exp \left( i \frac{1}{\sqrt{N}} \sum_{i=1}^N \left( X_i -m_i \right) \xi \right) \right) \\
&\qquad = \cov_{\mu ^{\sigma}} \left(  \mathbb{E}_{\mu_c ^{\sigma}} \left[ h(X) \right]  ,  \mathbb{E}_{\mu_c ^{\sigma}} \left[ \exp \left( i \frac{1}{\sqrt{N}} \sum_{i=1}^N \left( X_i -m_i \right) \xi \right) \right] \right) \label{e_total_covariance_first} \\
& \qquad \quad + \mathbb{E}_{\mu^{\sigma}} \left[ \cov_{\mu_c ^{\sigma}} \left( h(X) ,  \exp \left( i \frac{1}{\sqrt{N}} \sum_{i=1}^N \left( X_i -m_i \right) \xi \right) \right) \right]. \label{e_total_covariance_second}
\end{align}
We estimate~\eqref{e_total_covariance_first} and~\eqref{e_total_covariance_second} separately. \\

\begin{lemma} \label{l_main_equiv_ob_first}
Under the same assumptions as in Lemma~\ref{l_1st_der_inner}, it holds that
\begin{align}
\left| T_{\eqref{e_total_covariance_first}} \right|  \lesssim \| \nabla f \|_{\infty}|\supp f |^{\frac{1}{2}} \exp \left( -C M \ln N\right) \exp \left( - C \xi ^2 \right).
\end{align}
\end{lemma}

\medskip

\noindent \emph{Proof of Lemma~\ref{l_main_equiv_ob_first}.} \ Let us further decompose the set~$F_S$ into the boundary set (with respect to~$E_S$)~$F_S ^1$ and the exterior set~$F_S ^2$ as follows:
\begin{align} 
F_S ^1 : = \{ i \in F_S : \text{dist}(i, E_S) \leq R  \}, \label{e_boundary} \\
F_S ^2 : = \{ i \in F_S : \text{dist}(i, E_S) > R  \}. \label{e_exterior}
\end{align}
We note that~$[N]$ is decomposed into
\begin{align}
    [N] = E_s \cup F_S ^1 \cup F_S ^2,
\end{align}
where the union is disjoint. We also denote for each~$i \in E_S$
\begin{align} \label{e_conditional_mean_spin}
\tilde{m}_i : = \mathbb{E}_{\mu_c ^{\sigma}} \left[ X_i \right].
\end{align}
We write
\begin{align}
&\mathbb{E}_{\mu_c ^{\sigma}} \left[ \exp \left( i \frac{1}{\sqrt{N}} \sum_{i=1}^{N} \left( X_i -m_i \right) \xi \right) \right] \\
&\qquad = \exp \left( i \frac{1}{\sqrt{N}} \sum_{j \in F_S ^1} \left( X_j -m_j \right) \xi \right) \cdot \exp \left( i \frac{1}{\sqrt{N}} \sum_{j \in F_S ^2} \left( X_j -m_j \right) \xi \right) \\
& \qquad \qquad \times \exp \left( i \frac{1}{\sqrt{N}} \sum_{i \in E_S} \left( \tilde{m}_i -m_i \right) \xi \right) \cdot \mathbb{E}_{\mu_c ^{\sigma}} \left[ \exp \left( i \frac{1}{\sqrt{N}} \sum_{i \in E_S} \left( X_i -\tilde{m}_i \right) \xi \right) \right] \\
& \qquad = A \cdot B \cdot C \cdot D. \label{e_abcd}
\end{align}
Due to the finite range interaction (with interaction range~$R$), the conditional expectations
\begin{align}
\tilde{m}_i =\mathbb{E}_{\mu_c ^{\sigma}} \left[ X_i \right], \  \mathbb{E}_{\mu_c ^{\sigma}} \left[ h(X) \right], \text{ and } \mathbb{E}_{\mu_c ^{\sigma}} \left[ \exp \left( i \frac{1}{\sqrt{N}} \sum_{i \in E_S} \left( X_i -\tilde{m}_i \right) \xi \right) \right]
\end{align}
are only dependent on spins at~$F_S ^1$ (and thus independent of spins at~$F_S ^2$). In particular,~$A, C, D$ from~\eqref{e_abcd} and~$\mathbb{E}_{\mu_c ^{\sigma}} \left[ h(X) \right]$ are only dependent on spins at~$F_S ^1$. Thus we have
\begin{align}
T_{\eqref{e_total_covariance_first}} &= \mathbb{E}_{\mu^{\sigma}}\left[ \left(\mathbb{E}_{\mu_c ^{\sigma}} \left[ h(X) \right] - \mathbb{E}_{\mu^{\sigma}} \left[ h(X) \right]  \vphantom{\exp \left( i \frac{1}{\sqrt{N}} \sum_{j \in F_S ^1} \left( X_j-m_j \right)\xi + i \frac{1}{\sqrt{N}} \sum_{i \in E_S} \left( \tilde{m}_i - m_i \right) \xi \right)}  \right) A \cdot B \cdot C \cdot D  \right] \\
& = \mathbb{E}_{\mu^{\sigma}}\left[ \mathbb{E}_{\mu^{\sigma}} \left. \left[   \left(\mathbb{E}_{\mu_c ^{\sigma}} \left[ h(X) \right] - \mathbb{E}_{\mu^{\sigma}} \left[ h(X) \right] \vphantom{\exp \left( i \frac{1}{\sqrt{N}} \sum_{j \in F_S ^1} \left( X_j-m_j \right)\xi + i \frac{1}{\sqrt{N}} \sum_{i \in E_S} \left( \tilde{m}_i - m_i \right) \xi \right)} \right) A \cdot B \cdot C \cdot D  \ \right| \ X_i, i \in E_S \cup F_S ^1  \right] \  \right] \\
& = \mathbb{E}_{\mu^{\sigma}} \left[  \left(\mathbb{E}_{\mu_c ^{\sigma}} \left[ h(X) \right] - \mathbb{E}_{\mu^{\sigma}} \left[ h(X) \right] \vphantom{\exp \left( i \frac{1}{\sqrt{N}} \sum_{j \in F_S ^1} \left( X_j-m_j \right)\xi + i \frac{1}{\sqrt{N}} \sum_{i \in E_S} \left( \tilde{m}_i - m_i \right) \xi \right)} \right) A \cdot C \cdot D \cdot  \mathbb{E}_{\mu^{\sigma}} \left. \left[ B \ \right| \ X_i , i \in E_S \cup F_S ^1 \right] \  \right] \\
& = \mathbb{E}_{\mu^{\sigma}} \left[ \exp \left( i \frac{1}{\sqrt{N}} \sum_{j \in F_S ^1} \left( X_j-m_j \right)\xi + i \frac{1}{\sqrt{N}} \sum_{i \in E_S} \left( \tilde{m}_i - m_i \right) \xi \right)  \right. \\
& \qquad \qquad \quad  \times \mathbb{E}_{\mu_c ^{\sigma}} \left[ \exp \left( i \frac{1}{\sqrt{N}} \sum_{i \in E_S} \left( X_i -\tilde{m}_i \right) \xi \right) \right] \left(\mathbb{E}_{\mu_c ^{\sigma}} \left[ h(X) \right] - \mathbb{E}_{\mu^{\sigma}} \left[ h(X) \right] \vphantom{\exp \left( i \frac{1}{\sqrt{N}} \sum_{j \in F_S ^1} \left( X_j-m_j \right)\xi + i \frac{1}{\sqrt{N}} \sum_{i \in E_S} \left( \tilde{m}_i - m_i \right) \xi \right)} \right) \\
& \qquad \qquad \qquad \quad \quad  \times \left. \mathbb{E}_{\mu^{\sigma}} \left[ \left. \exp \left( i \frac{1}{\sqrt{N}} \sum_{k \in F_S ^2} \left( X_j-m_j \right)\xi \right) \ \right| \  X_i, i \in E_S \cup F_S ^1  \right] \ \right]. \label{e_total_cov_first_decomposition}
\end{align}
It holds by Lemma~\ref{l_conditional_function_difference} that
\begin{align}
&\left| \mathbb{E}_{\mu_c ^{\sigma}} \left[ h(X) \right] - \mathbb{E}_{\mu^{\sigma}} \left[ h(X) \right] \right|\\
& \qquad = \left| \int \left( \mathbb{E}_{\mu ^{\sigma} (dx^{E_S}  |  y^{F_S}) } \left[ h(X) \right] - \mathbb{E}_{\mu ^{\sigma} (dx^{E_S}  |  z^{F_S}) } \left[ h(X) \right] \right) \bar{\mu}^{\sigma} (dz^{F_S} ) \right|\\
& \qquad \lesssim \int  \| \nabla h \|_{\infty}  \left( \sum_{\substack{i \in E_S , j \in F_S \\ |i-j| \leq R}} M_{ij}^2(y_j -z_j )^2 \right)^{\frac{1}{2}} \exp \left( -C M \ln N\right) \bar{\mu}^{\sigma} (dz^{F_S} ). \label{e_application_conditional_function_diff}
\end{align}
Then a combination of~\eqref{e_total_cov_first_decomposition},~\eqref{e_application_conditional_function_diff} and Lemma~\ref{l_exchanging_exponential_terms} from below yields
\begin{align}
\left| T_{\eqref{e_total_covariance_first}} \right| & \lesssim \| \nabla h \|_{\infty} \exp \left( -C M \ln N\right) \exp \left( - C \xi ^2 \right) \\
& \qquad \times \mathbb{E}_{\mu^{\sigma}} \left[ \int   \left( \sum_{\substack{i \in E_S , j \in F_S \\ |i-j| \leq R}} M_{ij}^2(y_j -z_j )^2 \right)^{\frac{1}{2}}  \bar{\mu}^{\sigma} (dz^{F_S} ) \right]. \label{e_tot_cov_first_term_last_est}
\end{align}
Because there are at most~$2R^2 |\supp f | \sim |\supp f| $ many pairs of~$(i,j)$ with~$i \in E_S$, $j \in F_S$ with $ |i-j| \leq R$, an application of Schwarz inequality implies, as desired,
\begin{align}
&\mathbb{E}_{\mu^{\sigma}} \left[ \int   \left( \sum_{\substack{i \in E_S , j \in F_S \\ |i-j| \leq R}} M_{ij}^2(y_j -z_j )^2 \right)^{\frac{1}{2}}  \bar{\mu}^{\sigma} (dz^{F_S} ) \right] \\
& \qquad \lesssim \left( \int \int \sum_{\substack{i \in E_S , j \in F_S \\ |i-j| \leq R}} M_{ij}^2(y_j -z_j )^2   \bar{\mu}^{\sigma} (dz^{F_S} )  \mu^{\sigma} (dy^{F_S} )\right)^{\frac{1}{2}} \\
& \qquad =  \left( \sum_{\substack{i \in E_S , j \in F_S \\ |i-j| \leq R}} \int \int  M_{ij}^2(y_j -z_j )^2   \mu^{\sigma} (dz^{F_S} )  \mu^{\sigma} (dy^{F_S} )\right)^{\frac{1}{2}} \\
& \qquad \lesssim \left( \sum_{\substack{i \in E_S , j \in F_S \\ |i-j| \leq R}} 2 \var_{\mu^{\sigma}} \left( X_j \right)   \right) ^{\frac{1}{2}} \overset{Corollary~\ref{l_moment_estimate}}{\lesssim} |\supp f|^{\frac{1}{2}}. \label{e_mij_proof}
\end{align}
By the definition~\eqref{e_assumption_f_structure} and the inequality~\eqref{e_alpha_g_estimate}, we have
\begin{align}
\| \nabla h \|_{\infty} \leq \| \nabla f  \|_{\infty} + |c_f| \left\| \nabla \left( \sum_{i \in S} X_i \right) \right\|_{\infty} \lesssim \| \nabla f \|_{\infty}. \label{e_comparison_fh_gradient}
\end{align}
It also holds from definition~\eqref{e_assumption_f_structure} that~$|\supp f | = |\supp h|$. Therefore we conclude from~\eqref{e_tot_cov_first_term_last_est},~\eqref{e_mij_proof} and~\eqref{e_comparison_fh_gradient} that
\begin{align}
\left| T_{\eqref{e_total_covariance_first}} \right| \lesssim \| \nabla f \|_{\infty}|\supp f |^{\frac{1}{2}} \exp \left( -C M \ln N\right) \exp \left( - C \xi ^2 \right).
\end{align}
\qed
\medskip

To estimate~\eqref{e_total_covariance_second}, we need the following extension of~\cite[Lemma 7]{KwMe18} (see also~\cite[Lemma 10]{kwme18a}).

\begin{lemma}[Extension of Lemma 7 in~\cite{KwMe18}] \label{l_exchanging_exponential_terms}
For large enough~$N$ and~$\delta>0$ small enough, there exists a positive constant~$C>0$ such that the following inequalities hold for all~$\xi \in \mathbb{R}$ with~$ \frac{\left|\xi\right|}{\sqrt{N}}  \leq \delta$.
\begin{align}
 \left|\mathbb{E}_{\mu^{\sigma}} \left[ \left. \exp \left( i \frac{1}{\sqrt{N}} \sum_{k \in F_S ^2} \left( X_j-m_j \right)\xi \right) \ \right| \  X_i, i \in E_S \cup F_S ^1  \right] \right| \lesssim \exp\left(-C\xi^2 \right) .
\label{e_lemma exchanging 2nd der}
\end{align}
\end{lemma}

\medskip

\begin{remark}
The proof of Lemma~\ref{l_exchanging_exponential_terms} is almost similar to that of~\cite[Lemma 7]{KwMe18}. One should compare the sets~$\left( E_S \cup F_S ^1, F_S ^2 \right)$ with~$( F_1^{n,l} , F_2^{n,l} )$ in~\cite{KwMe18}. The main difference is that we assume finite range interaction with range~$R$ instead of the nearest neighbor interaction. However, there is only a cosmetic difference between these two proofs. We leave the details to the reader. 
\end{remark}

\medskip

The next statement is an estimation of~\eqref{e_total_covariance_second}. \\

\begin{lemma} \label{l_main_equiv_ob_second_aux}
Under the same settings as in Lemma~\ref{l_1st_der_inner}, it holds that
\end{lemma}
\begin{align}
\left| T_{\eqref{e_total_covariance_second}} \right| & \lesssim \| \nabla f \|_{\infty} \left(1+ \int \left( \sum_{\substack{i \in E_S , j \in F_S \\ |i-j| \leq R}} M_{ij}^2(y_j -z_j )^2 \right)^{\frac{1}{2}} \bar{\mu}^{\sigma} (dz) \right)\frac{|\xi|}{N} \\
& \quad + \| \nabla f \|_{\infty} \frac{|\supp f |}{N} \xi^2  + \| \nabla f \|_{\infty} \frac{ |\supp f |}{N} |\xi|^3.
\end{align}

\medskip

\noindent \emph{Proof of Lemma~\ref{l_main_equiv_ob_second_aux}.} \ Let us recall the definition~\eqref{e_boundary} and~\eqref{e_exterior} of~$F_S ^1$ and~$F_S ^2$, respectively. As in Lemma~\ref{l_main_equiv_ob_first}, we write~\eqref{e_total_covariance_second} as (see~\eqref{e_total_cov_first_decomposition})
\begin{align}
T_{\eqref{e_total_covariance_second}} & = \mathbb{E}_{\mu^{\sigma}} \left[ \exp \left( i \frac{1}{\sqrt{N}} \sum_{j \in F_S ^1} \left( X_j-m_j \right)\xi + i \frac{1}{\sqrt{N}} \sum_{i \in E_S} \left( \tilde{m}_i - m_i \right) \xi \right)  \right. \\
& \qquad \qquad \quad \times \cov_{\mu_c ^{\sigma}} \left( h(X)  ,   \exp \left( i \frac{1}{\sqrt{N}} \sum_{i \in E_S} \left( X_i -\tilde{m}_i \right) \xi \right) \right) \\
& \qquad \qquad \qquad \quad \quad \times \left. \mathbb{E}_{\mu^{\sigma}} \left[ \left. \exp \left( i \frac{1}{\sqrt{N}} \sum_{k \in F_S ^2} \left( X_j-m_j \right)\xi \right) \ \right| \  X_i, i \in E_S \cup F_S ^1  \right] \right]. \label{e_total_cov_second_decomposition}
\end{align}
We then apply third order Taylor expansions to get
\begin{align}
&\cov_{\mu_c ^{\sigma}} \left( h(X)  ,   \exp \left( i \frac{1}{\sqrt{N}} \sum_{i \in E_S} \left( X_i -\tilde{m}_i \right) \xi \right) \right) \label{e_original} \\
& = \cov_{\mu_c ^{\sigma}} \left( h(X)  ,   \sum_{i \in E_S} \left( X_i -\tilde{m}_i \right)  \right) i \frac{1}{\sqrt{N}} \xi \label{e_linear}\\
&  \quad + \frac{1}{2} \cov_{\mu_c ^{\sigma}} \left( h(X)  ,   \left(\sum_{i \in E_S} \left( X_i -\tilde{m}_i \right)\right)^2  \right) \left(i \frac{1}{\sqrt{N}} \xi \right)^2 \label{e_quadratic} \\
&  \quad + \frac{1}{6} \cov_{\mu_c ^{\sigma}} \left( h(X)  ,   \left(\sum_{i \in E_S} \left( X_i -\tilde{m}_i \right)\right)^3 \exp \left( i \frac{1}{\sqrt{N}} \sum_{i \in E_S} \left( X_i -\tilde{m}_i \right) \tilde{\xi} \right)  \right) \left(i \frac{1}{\sqrt{N}} \xi \right)^3, \label{e_cubic}
\end{align}
where~$\tilde{\xi}$ is a real number between~$0$ and~$\xi$. Let us begin with estimation of~\eqref{e_linear}. Recalling the definition~\eqref{e_assumption_f_structure} of the function~$h$, it holds that
\begin{align}
\left|T_{\eqref{e_linear}}\right| &= \left|\cov_{\mu_c ^{\sigma}} \left( f  ,   \sum_{i \in E_S} X_i  \right) - c_f \cov_{\mu_c ^{\sigma}} \left( \sum_{i \in S}\sum_{k \in B_{\varphi(i)}} X_k  ,   \sum_{i \in E_S} X_i \right) \right| \frac{|\xi|}{\sqrt{N}} \\
& \leq \left|  \cov_{\mu_c ^{\sigma}} \left( f  ,   \sum_{i \in E_S} X_i  \right) - \cov_{\mu ^{\sigma}} \left( f  ,   \sum_{i \in E_S} X_i  \right) \right| \frac{|\xi|}{\sqrt{N}} \\
& \quad + \left| \cov_{\mu ^{\sigma}} \left( f  ,   \sum_{i \in E_S} X_i  \right)- c_f \cov_{\mu_c ^{\sigma}} \left( \sum_{i \in S}\sum_{k \in B_{\varphi(i)}} X_k  ,   \sum_{i \in E_S} X_i \right) \right| \frac{|\xi|}{\sqrt{N}} \\
&=  \left|  \cov_{\mu_c ^{\sigma}} \left( f  ,   \sum_{i \in E_S} X_i  \right) - \cov_{\mu ^{\sigma}} \left( f  ,   \sum_{i \in E_S} X_i  \right) \right| \frac{|\xi|}{\sqrt{N}} \label{e_linear_1} \\
& \quad + \left|c_f\right| \left|  \cov_{\mu ^{\sigma}} \left( \sum_{i \in S}\sum_{k \in B_{\varphi(i)}} X_k  ,   \sum_{j \in E_S} X_j \right) - \cov_{\mu_c ^{\sigma}} \left( \sum_{i \in S}\sum_{k \in B_{\varphi(i)}} X_k  ,   \sum_{j \in E_S} X_j \right) \right| \frac{|\xi|}{\sqrt{N}}. \label{e_linear_2}
\end{align}
Corollary~\ref{l_conditional_covariance_corollary} implies that
\begin{align}
T_{\eqref{e_linear_1}} &\lesssim \| \nabla f \|_{\infty} \exp \left( -C M \ln N \right) \frac{|\xi|}{\sqrt{N}} \\
& \quad + \left|E_S \right| \| \nabla f \|_{\infty} \exp \left( -C M \ln N\right) \int \left( \sum_{\substack{i \in E_S , j \in F_S \\ |i-j| \leq R}} M_{ij}^2(y_j -z_j )^2 \right)^{\frac{1}{2}} \bar{\mu}^{\sigma} (dz) \frac{|\xi|}{\sqrt{N}}\\
&  \quad + \left|E_S \right| \| \nabla f \|_{\infty}  |\supp f |^{\frac{1}{2}} \exp \left( - C M \ln N\right)\frac{|\xi|}{\sqrt{N}}.
\end{align}
Because~$|E_S | \leq 2 |\supp f | M \ln N$, it holds for~$N$ large enough that
\begin{align}
\left|T_{\eqref{e_linear_1}}\right| \lesssim \| \nabla f \|_{\infty} \left(1+ \int \left( \sum_{\substack{i \in E_S , j \in F_S \\ |i-j| \leq R}} M_{ij}^2(y_j -z_j )^2 \right)^{\frac{1}{2}} \bar{\mu}^{\sigma} (dz) \right) \frac{|\xi|}{N}. 
\end{align}
Similarly, using Corollary~\ref{l_conditional_covariance_corollary} and~\eqref{e_alpha_g_estimate}, we get
\begin{align}
\left|T_{\eqref{e_linear_2}} \right| &\lesssim \| \nabla f \|_{\infty} \left(1+ \int \left( \sum_{\substack{i \in E_S , j \in F_S \\ |i-j| \leq R}} M_{ij}^2(y_j -z_j )^2 \right)^{\frac{1}{2}} \bar{\mu}^{\sigma} (dz) \right) \frac{|\xi|}{N}. 
\end{align}
Therefore
\begin{align}
\left|T_{\eqref{e_linear}}\right| &\leq \left|T_{\eqref{e_linear_1}} \right| + \left|T_{\eqref{e_linear_2}} \right| \\
&\lesssim  \| \nabla f \|_{\infty} \left(1+ \int \left( \sum_{\substack{i \in E_S , j \in F_S \\ |i-j| \leq R}} M_{ij}^2(y_j -z_j )^2 \right)^{\frac{1}{2}} \bar{\mu}^{\sigma} (dz) \right) \frac{|\xi|}{N}.
\end{align}

The estimate for~\eqref{e_quadratic} follows from Lemma~\ref{l_function_second_estimate}:
\begin{align}
\left|T_{\eqref{e_quadratic}} \right| \lesssim \| \nabla h \|_{L^4 (\mu_c ^{\sigma})} |\supp h | \frac{\xi^2}{N} \leq \| \nabla h \|_{\infty} \frac{|\supp h |}{N} \xi^2 \overset{\eqref{e_comparison_fh_gradient}}{\lesssim} \| \nabla f \|_{\infty} \frac{|\supp f |}{N} \xi^2.
\end{align}

Let us turn to the estimation of~\eqref{e_cubic}. By applying H\"older's inequality we have
\begin{align}
&\left| \cov_{\mu_c ^{\sigma}} \left( h(X)  ,   \left(\sum_{i \in E_S} \left( X_i -\tilde{m}_i \right)\right)^3 \exp \left( i \frac{1}{\sqrt{N}} \sum_{i \in E_S} \left( X_i -\tilde{m}_i \right) \tilde{\xi} \right)  \right) \right|\\
& = \left|\mathbb{E}_{\mu_c ^{\sigma}} \left[ \left(  h(X) - \mathbb{E}_{\mu_c ^{\sigma}} \left[ h(X) \right]\right)\left(\sum_{i \in E_S} \left( X_i -\tilde{m}_i \right)\right)^3 \exp \left( i \frac{1}{\sqrt{N}} \sum_{i \in E_S} \left( X_i -\tilde{m}_i \right) \tilde{\xi} \right)   \right]\right| \\
& \leq \| h(X) -  \mathbb{E}_{\mu_c^{\sigma}} \left[ h(X) \right] \|_{L^4 ( \mu_c ^{\sigma})} \left(\mathbb{E}_{\mu_c ^{\sigma}}\left[\left(\sum_{i \in E_S} \left( X_i -\tilde{m}_i \right)\right)^4 \right] \right)^{\frac{3}{4}}. \label{e_holder}
\end{align}
A combination of~\eqref{e_alpha_g_estimate} and Lemma~\ref{l_fourth_moment} yields
\begin{align}
&\| h(X) - \mathbb{E}_{\mu_c ^{\sigma}} \left[ h(X) \right] \|_{L^4 (\mu_c ^{\sigma})} \\
&\qquad \leq \| f(X) - \mathbb{E}_{\mu_c ^{\sigma}} \left[ f(X) \right] \|_{L^4 (\mu_c ^{\sigma})} + |c_f| \left\| \sum_{i \in S}\sum_{k \in B_{\varphi(i)}} \left(X_k - \tilde{m}_k \right)\right\|_{L^4 (\mu_c ^{\sigma})} \\
& \qquad \lesssim  \|\nabla f \|_{\infty} + \frac{\| \nabla f \|_{\infty}}{|\supp f |^{\frac{1}{2}}} |\supp f |^{\frac{1}{2}} \lesssim \| \nabla f \|_{\infty} . \label{e_comparison_fh_fourth_moment}
\end{align}
Thus we conclude from~\eqref{e_holder},~\eqref{e_comparison_fh_fourth_moment} and Lemma~\ref{l_fourth_moment} that
\begin{align}
\left|T_{\eqref{e_cubic}} \right|
 \overset{\eqref{e_comparison_fh_fourth_moment}, \ Lemma~\ref{l_fourth_moment}}{\lesssim} \| \nabla f \|_{\infty} |E_S|^{\frac{3}{2}} \frac{ |\xi|^3}{N^{\frac{3}{2}}} \lesssim \| \nabla f \|_{\infty} \frac{ |\supp f|}{N} |\xi| ^3,
\end{align}
where we used~$|E_S| \leq 2|\supp f | M \ln N$ and thus for~$N$ large,
\begin{align}
\frac{|E_S|^{\frac{3}{2}}}{N^{\frac{3}{2}}} \lesssim \frac{|\supp f |}{N}.
\end{align}
\medskip 

Collecting all the estimates we have proven so far, we get
\begin{align}
\left| T_{\eqref{e_original}} \right| & \leq \left| T_{\eqref{e_linear}}\right| + \left| T_{\eqref{e_quadratic}}\right|+\left| T_{\eqref{e_cubic}}\right| \\
& \lesssim \| \nabla f \|_{\infty} \left(1+ \int \left( \sum_{\substack{i \in E_S , j \in F_S \\ |i-j| \leq R}} M_{ij}^2(y_j -z_j )^2 \right)^{\frac{1}{2}} \bar{\mu}^{\sigma} (dz) \right)\frac{|\xi|}{N} \\
& \quad + \| \nabla f \|_{\infty} \frac{|\supp f |}{N} \xi^2  + \| \nabla f \|_{\infty} \frac{ |\supp f |}{N} |\xi|^3.
\end{align}

\qed

\medskip 

Lemma~\ref{l_main_equiv_ob_second} is a direct consequence of Lemma~\ref{l_main_equiv_ob_second_aux}.
\begin{lemma} \label{l_main_equiv_ob_second}
Under the same settings as in Lemma~\ref{l_1st_der_inner}, it holds that
\begin{align}
&\left| \mathbb{E}_{\mu^{\sigma}} \left[ \cov_{\mu_c ^{\sigma}} \left( h(X)  ,  \exp \left( i \frac{1}{\sqrt{N}} \sum_{i=1}^N \left( X_i -m_i \right) \xi \right) \right) \right] \right| \lesssim \| \nabla f \|_{\infty} \frac{ |\supp f |}{N} \exp \left(- C\xi^2 \right).
\end{align}
\end{lemma}

\medskip

\noindent \emph{Proof of Lemma~\ref{l_main_equiv_ob_second}.} \ Let us recall the decomposition~\eqref{e_total_cov_second_decomposition}. We recall also estimation~\eqref{e_mij_proof}, which implies
\begin{align}
&\mathbb{E}_{\mu^{\sigma}} \left[ \int \left( \sum_{\substack{i \in E_S , j \in F_S \\ |i-j| \leq R}} M_{ij}^2(y_j -z_j )^2 \right)^{\frac{1}{2}} \bar{\mu}^{\sigma} (dz) \right] \lesssim |\supp f |^{\frac{1}{2}}, \label{e_mij_ineq}.
\end{align}
A combination of Lemma~\ref{l_exchanging_exponential_terms} and Lemma~\ref{l_main_equiv_ob_second_aux} yields
\begin{align}
&\left| \mathbb{E}_{\mu^{\sigma}} \left[ \cov_{\mu_c ^{\sigma}} \left( h(X)  ,  \exp \left( i \frac{1}{\sqrt{N}} \sum_{i=1}^N \left( X_i -m_i \right) \xi \right) \right) \right] \right| \\
& \qquad \lesssim \|\nabla f \|_{\infty} \frac{|\supp f |}{N}\left( |\xi| + \xi^2 + |\xi|^3 \right) \exp \left( - C\xi^2 \right)  \lesssim  \|\nabla f \|_{\infty} \frac{|\supp f |}{N} \exp \left( - C\xi^2 \right).
\end{align}
\qed

\medskip

Now we are ready to give a proof of Lemma~\ref{l_1st_der_inner}. \\

\noindent \emph{Proof of Lemma~\ref{l_1st_der_inner}.} \ The law of total covariance implies
\begin{align}
&\cov_{\mu^{\sigma}} \left( h(X)  ,   \exp \left( i \frac{1}{\sqrt{N}} \sum_{i=1}^N \left( X_i -m_i \right) \xi \right) \right) \\
&\qquad = \cov_{\mu ^{\sigma}} \left(  \mathbb{E}_{\mu_c ^{\sigma}} \left[ h(X) \right]  ,  \mathbb{E}_{\mu_c ^{\sigma}} \left[ \exp \left( i \frac{1}{\sqrt{N}} \sum_{i=1}^N \left( X_i -m_i \right) \xi \right) \right] \right) \\
& \qquad \quad + \mathbb{E}_{\mu^{\sigma}} \left[ \cov_{\mu_c ^{\sigma}} \left( h(X)  ,  \exp \left( i \frac{1}{\sqrt{N}} \sum_{i=1}^N \left( X_i -m_i \right) \xi \right) \right) \right]. 
\end{align}
By Lemma~\ref{l_main_equiv_ob_first} and Lemma~\ref{l_main_equiv_ob_second} it holds that for~$M$,~$N$ large enough,
\begin{align}
&\left| \int_{\{ \left| \left(1/ \sqrt{N}\right)\xi\right|\leq\delta\}} \cov_{\mu^{\sigma}} \left( h(X)  ,   \exp \left( i \frac{1}{\sqrt{N}} \sum_{i=1}^N \left( X_i -m_i \right) \xi \right) \right) d\xi \right| \\
& \qquad \lesssim \int_{\{ \left| \left(1/ \sqrt{N}\right)\xi\right|\leq\delta\}} \| \nabla f \|_{\infty}|\supp f |^{\frac{1}{2}} \exp \left( -C M \ln N\right) \exp \left( - C \xi ^2 \right) d\xi \\
& \qquad  \quad +  \int_{\{ \left| \left(1/ \sqrt{N}\right)\xi\right|\leq\delta\}} \| \nabla f \|_{\infty} \frac{ |\supp f |}{N} \exp \left(- C\xi^2 \right) d\xi \\
& \qquad \leq \int \| \nabla f \|_{\infty}|\supp f |^{\frac{1}{2}} \exp \left( -C M \ln N\right) \exp \left( - C \xi ^2 \right) d\xi \\
& \qquad \quad + \int \| \nabla f \|_{\infty} \frac{ |\supp f |}{N} \exp \left(- C\xi^2 \right) d\xi \\
& \qquad \lesssim \| \nabla f \|_{\infty}\frac{ |\supp f |}{N}.
\end{align}
\qed

\medskip 

\begin{remark} \label{r_hij}
A detailed review show that the arguments in this section can be adapted to yield similar results applied to~$h_{ij}$ (see~\eqref{e_def_hij}) instead of the function~$h$. The only place where one should check details is the proof of Lemma~\ref{l_main_equiv_ob_second_aux}, especially the estimation of~$T_{\eqref{e_linear}}$. We choose not to outline the details because they would yield many redundancies. 
\end{remark}
\medskip

\section{Proof of Proposition~\ref{p_2nd_der_main_proposition}} \label{s_proof_2nd_der_main_proposition}

Proof of Proposition~\ref{p_2nd_der_main_proposition} follows the same idea of Proposition~\ref{p_1st_der_main_proposition} with more careful estimation. We follow similar calculations as in the proof of~\cite[Proposition 3]{kwme18a}. Instead of a second order Taylor expansion we use this time a third order Taylor expansion, which leads to improved estimates.\\

In this section, the set~$S$ denotes union of~$\supp f$ and~$\supp g$, i.e.~$S = \supp f \cup \supp g$. Let us recall the definition~\eqref{d_def_e1f} and~\eqref{d_def_e2f} of the sets~$E_S$ and~$F_S$ and the decomposition~\eqref{e_measure_decomposition} of the gce~$\mu^{\sigma}$. We write
\begin{align}
\mu_c ^{\sigma} (dx| y ) = \mu^{\sigma}( dx^{E_S} | y^{F_S}), \qquad \bar{\mu}^{\sigma} (dy ) = \bar{\mu}^{\sigma} ( dy^{F_S}).
\end{align}

As before, the integral is divided into inner and outer parts and estimated separately. More precisely, let us fix~$\delta>0$ small enough and decompose the integral as
\begin{align}
& \int_{\mathbb{R}} \mathbb{E}_{\mu^{\sigma}} \left[ \left( f(X) - \mathbb{E}_{\mu^{\sigma}} \left[f(X) \right] \right)\left( g(X) - \mathbb{E}_{\mu^{\sigma}} \left[g(X) \right] \right)  \exp \left( i \frac{1}{\sqrt{N}} \sum_{k=1}^{N} \left( X_k - m_k \right) \xi \right) \right] d\xi \\
& \qquad = \int_{\{  \left| \left(1/\sqrt{N}\right) \xi \right| \leq \delta \}} \mathbb{E}_{\mu^{\sigma}} \left[ \exp \left( i \frac{1}{\sqrt{N}} \sum_{k=1}^{N} \left( X_k - m_k \right) \xi \right) \right. \\
&\qquad \qquad \qquad \qquad \qquad \qquad \qquad \left. \times \vphantom{\exp \left( i \frac{1}{\sqrt{N}} \sum_{k=1}^{N} \left( X_k - m_k \right) \xi \right)} \left( f(X) - \mathbb{E}_{\mu^{\sigma}} \left[f(X) \right] \right)\left( g(X) - \mathbb{E}_{\mu^{\sigma}} \left[g(X) \right] \right) \right] d\xi \label{e_exp_decay_inner2}\\
& \qquad \quad + \int_{\{  \left| \left(1/\sqrt{N}\right) \xi \right| > \delta \}} \mathbb{E}_{\mu^{\sigma}} \left[ \exp \left( i \frac{1}{\sqrt{N}} \sum_{k=1}^{N} \left( X_k - m_k \right) \xi \right) \right. \\
&\qquad \qquad \qquad \qquad \qquad \qquad \qquad \quad  \left. \times \vphantom{\exp \left( i \frac{1}{\sqrt{N}} \sum_{k=1}^{N} \left( X_k - m_k \right) \xi \right)} \left( f(X) - \mathbb{E}_{\mu^{\sigma}} \left[f(X) \right] \right)\left( g(X) - \mathbb{E}_{\mu^{\sigma}} \left[g(X) \right] \right) \right] d\xi \label{e_exp_decay_outer2}
\end{align}
As in the proof of Proposition~\ref{p_1st_der_main_proposition}, the estimation of outer integral~\eqref{e_exp_decay_outer2} follows from a slight modification of arguments in~\cite{kwme18a}. More precisely, we have
\begin{align} \label{e_outer_estimate2}
\left| T_{\eqref{e_exp_decay_outer2}} \right| \lesssim \| \nabla f \|_{\infty} \| \nabla g \|_{\infty} \left( \frac{|\supp f | + |\supp g | }{N}+  \exp\left(-Cd_{f,g}\right)  \right).
\end{align}

Let us state and prove the following lemma, which corresponds to Lemma~\ref{l_1st_der_inner} in the proof of Proposition~\ref{p_equivalence_decay_of_correlations}. \\

\begin{lemma} \label{l_2nd_der_inner} 
It holds that
\begin{align}
&\left| \int_{\{  \left| \left(1/\sqrt{N}\right) \xi \right| \leq \delta \}} \mathbb{E}_{\mu^{\sigma}} \left[ \exp \left( i \frac{1}{\sqrt{N}} \sum_{k=1}^{N} \left( X_k - m_k \right) \xi \right) \right. \right. \\
& \qquad \qquad \qquad \qquad \qquad \qquad \left. \left. \times \vphantom{\exp \left( i \frac{1}{\sqrt{N}} \sum_{k=1}^{N} \left( X_k - m_k \right) \xi \right)} \left( f(X) - \mathbb{E}_{\mu^{\sigma}} \left[f(X) \right] \right)\left( g(X) - \mathbb{E}_{\mu^{\sigma}} \left[g(X) \right] \right) \right] d\xi \right| \\
& \qquad \lesssim \| \nabla f \|_{\infty} \| \nabla g \|_{\infty} \left( \frac{|\supp f | + |\supp g | }{N}+  \exp\left(-Cd_{f,g}\right)  \right).
\end{align}
\end{lemma}

\medskip

\noindent \emph{Proof of Proposition~\ref{p_2nd_der_main_proposition}.} \ This directly follows from Lemma~\ref{l_2nd_der_inner} and~\eqref{e_outer_estimate2}. \qed  

\medskip

To prove Lemma~\ref{l_2nd_der_inner}, we first write
\begin{align}
&\mathbb{E}_{\mu^{\sigma}  } \left[ \exp\left( i\frac{1}{\sqrt{N}} \sum_{k=1}^N \left(X_k -m_k \right) \xi \right)\left(f(X)-\mathbb{E}_{\mu^{\sigma}} \left[f(X) \right] \right)\left(g(X)-\mathbb{E}_{\mu^{\sigma}} \left[g(X) \right] \right) \right] \\
&\qquad  = \mathbb{E}_{\mu^{\sigma}  } \left[ \exp\left( i\frac{1}{\sqrt{N}} \sum_{k=1}^N \left(X_k -m_k \right) \xi \right) \right. \\
&\left. \qquad \qquad \qquad \qquad \quad \times  \vphantom{\exp\left( i\frac{1}{\sqrt{N}} \sum_{k=1}^N \left(X_k -m_k \right) \xi \right)}  \left(f(X)-\mathbb{E}_{\mu_c^{\sigma}} \left[f(X) \right] \right)\left(g(X)-\mathbb{E}_{\mu_c^{\sigma}} \left[g(X) \right] \right) \right] \label{e_cov_decomp_1} \\
&\qquad \quad  + \mathbb{E}_{\mu^{\sigma}  } \left[ \exp\left( i\frac{1}{\sqrt{N}} \sum_{k=1}^N \left(X_k -m_k \right) \xi \right) \right. \\
& \left. \qquad \qquad \qquad \qquad \quad \times \vphantom{\exp\left( i\frac{1}{\sqrt{N}} \sum_{k=1}^N \left(X_k -m_k \right) \xi \right)} \left(\mathbb{E}_{\mu_c^{\sigma}} \left[f(X) \right]-\mathbb{E}_{\mu^{\sigma}} \left[f(X) \right] \right)\left(g(X)-\mathbb{E}_{\mu_c^{\sigma}} \left[g(X) \right] \right) \right] \label{e_cov_decomp_2}\\
& \qquad \quad + \mathbb{E}_{\mu^{\sigma}  } \left[ \exp\left( i\frac{1}{\sqrt{N}} \sum_{k=1}^N \left(X_k -m_k \right) \xi \right) \right. \\
& \left. \qquad \qquad \qquad \qquad \quad \times \vphantom{\exp\left( i\frac{1}{\sqrt{N}} \sum_{k=1}^N \left(X_k -m_k \right) \xi \right)}  \left(f(X) -\mathbb{E}_{\mu_c^{\sigma}} \left[f(X) \right] \right)\left(\mathbb{E}_{\mu_c^{\sigma}} \left[g(X) \right]-\mathbb{E}_{\mu^{\sigma}} \left[g(X) \right] \right) \right] \label{e_cov_decomp_3}\\
& \qquad \quad + \mathbb{E}_{\mu^{\sigma}  } \left[ \exp\left( i\frac{1}{\sqrt{N}} \sum_{k=1}^N \left(X_k -m_k \right) \xi \right) \right. \\
&\left. \qquad \qquad \qquad \qquad \quad \times  \vphantom{\exp\left( i\frac{1}{\sqrt{N}} \sum_{k=1}^N \left(X_k -m_k \right) \xi \right)}  \left(\mathbb{E}_{\mu_c^{\sigma}} \left[f(X) \right]-\mathbb{E}_{\mu^{\sigma}} \left[f(X) \right] \right)\left(\mathbb{E}_{\mu_c^{\sigma}} \left[g(X) \right]-\mathbb{E}_{\mu^{\sigma}} \left[g(X) \right] \right) \right]. \label{e_cov_decomp_4}
\end{align}

\medskip 

\begin{lemma} \label{l_2nd_der_inner_1}
It holds that
\begin{align}
\left|T_{\eqref{e_cov_decomp_1}} \right| \lesssim \| \nabla f \|_{\infty} \| \nabla g \|_{\infty} \left( \frac{|\supp f | + |\supp g | }{N}+  \exp\left(-Cd_{f,g}\right)  \right) \exp \left( -C \xi^2 \right).
\end{align}
\end{lemma}

\medskip

\noindent \emph{Proof of Lemma~\ref{l_2nd_der_inner_1}} \ Let us recall the decomposition~\eqref{e_boundary},~\eqref{e_exterior} of~$F_S^1$,~$F_S ^2$ and the definition~\eqref{e_conditional_mean_spin} of~$\tilde{m}_i$. It holds that (see for example~\eqref{e_total_cov_first_decomposition} and~\eqref{e_total_cov_second_decomposition})
\begin{align}
T_{\eqref{e_cov_decomp_1}}& = \mathbb{E}_{\mu^{\sigma}} \left[ \exp \left( i \frac{1}{\sqrt{N}} \sum_{j \in F_S ^1} \left( X_j-m_j \right)\xi + i \frac{1}{\sqrt{N}} \sum_{i \in E_S} \left( \tilde{m}_i - m_i \right) \xi \right)  \right. \\
& \qquad \qquad \quad  \times \mathbb{E}_{\mu_c ^{\sigma}} \left[ \exp \left( i \frac{1}{\sqrt{N}} \sum_{i \in E_S} \left( X_i -\tilde{m}_i \right) \xi \right)   \right. \\
&\left. \qquad \qquad \qquad \qquad \qquad  \times \vphantom{\exp\left( i\frac{1}{\sqrt{N}} \sum_{k=1}^N \left(X_k -m_k \right) \xi \right)} \left(f(X)-\mathbb{E}_{\mu_c^{\sigma}} \left[f(X) \right] \right)\left(g(X)-\mathbb{E}_{\mu_c^{\sigma}} \left[g(X) \right] \right)   \right] \\
& \qquad \qquad \qquad \quad \quad  \times \left. \mathbb{E}_{\mu^{\sigma}} \left[ \left. \exp \left( i \frac{1}{\sqrt{N}} \sum_{k \in F_S ^2} \left( X_j-m_j \right)\xi \right)  \right| \  X_i, i \in E_S \cup F_S ^1  \right]  \right], \label{e_second_der_further_decomposition_1}
\end{align}
Taylor expansion implies that there is~$\tilde{\xi}$ between~$0$ and~$\xi$ such that
\begin{align}
&\mathbb{E}_{\mu_c ^{\sigma}} \left[\left(f(X)-\mathbb{E}_{\mu_c^{\sigma}} \left[f(X) \right] \right)\left(g(X)-\mathbb{E}_{\mu_c^{\sigma}} \left[g(X) \right] \right)\exp \left( i \frac{1}{\sqrt{N}} \sum_{i \in E_S} \left( X_i -\tilde{m}_i \right) \xi \right)  \right] \label{e_second_taylor_original} \\
& = \mathbb{E}_{\mu_c ^{\sigma}} \left[\left(f(X)-\mathbb{E}_{\mu_c^{\sigma}} \left[f(X) \right] \right)\left(g(X)-\mathbb{E}_{\mu_c^{\sigma}} \left[g(X) \right] \right)  \right] \label{e_second_taylor_0}  \\
& \quad + \mathbb{E}_{\mu_c ^{\sigma}} \left[\left(f(X)-\mathbb{E}_{\mu_c^{\sigma}} \left[f(X) \right] \right)\left(g(X)-\mathbb{E}_{\mu_c^{\sigma}} \left[g(X) \right] \right)\sum_{i \in E_S} \left( X_i -\tilde{m}_i \right)   \right] \left(i \frac{\xi}{\sqrt{N}} \right) \label{e_second_taylor_1} \\
& \quad + \frac{1}{2} \mathbb{E}_{\mu_c ^{\sigma}} \left[\left(f(X)-\mathbb{E}_{\mu_c^{\sigma}} \left[f(X) \right] \right)\left(g(X)-\mathbb{E}_{\mu_c^{\sigma}} \left[g(X) \right] \right) \left(\sum_{i \in E_S} \left( X_i -\tilde{m}_i \right) \right)^2  \right] \left(i \frac{\xi}{\sqrt{N}} \right)^2 \label{e_second_taylor_2}\\
&  \quad + \frac{1}{6} \mathbb{E}_{\mu_c ^{\sigma}} \left[\left(f(X)-\mathbb{E}_{\mu_c^{\sigma}} \left[f(X) \right] \right)\left(g(X)-\mathbb{E}_{\mu_c^{\sigma}} \left[g(X) \right] \right)  \vphantom{\left(\sum_{i \in E_S} \left( X_i -\tilde{m}_i \right) \right)^3}  \right. \\
&\left. \qquad \qquad \qquad \quad \times  \left(\sum_{i \in E_S} \left( X_i -\tilde{m}_i \right) \right)^3 \exp \left( i \frac{1}{\sqrt{N}} \sum_{i \in E_S} \left( X_i -\tilde{m}_i \right) \tilde{\xi} \right)  \right] \left(i \frac{\xi}{\sqrt{N}} \right)^3.  \label{e_second_taylor_3}
\end{align}
Let us begin with estimation of~\eqref{e_second_taylor_0}. An application of Theorem~\ref{p_decay_of_correlations_gce} implies that
\begin{align}
\left| T_{\eqref{e_second_taylor_0}}\right| \lesssim \|\nabla f \|_{L^2 (\mu_c ^{\sigma})}\|\nabla g \|_{L^2 (\mu_c ^{\sigma})} \exp \left( - C d_{f,g} \right) \leq \|\nabla f \|_{\infty}\|\nabla g \|_{\infty} \exp \left( - C d_{f,g} \right).
\end{align}
Next, to estimate~\eqref{e_second_taylor_1}, let us decompose~$E_S$ into two parts:
\begin{align}
E_S ^f : = \{ i \in E_S \ : \ \text{dist}(i \ , \ \supp f ) \geq \frac{1}{2} d_{f,g} \},
\qquad E_S ^g : = E_S \backslash E_S ^f.
\end{align}
We observe that for each~$i \in E_S ^g$,
\begin{align}
\text{dist}(i, \supp g ) \geq \frac{1}{2} d_{f,g} .
\end{align}
We write
\begin{align}
&\mathbb{E}_{\mu_c ^{\sigma}} \left[\left(f(X)-\mathbb{E}_{\mu_c^{\sigma}} \left[f(X) \right] \right)\left(g(X)-\mathbb{E}_{\mu_c^{\sigma}} \left[g(X) \right] \right) \sum_{i \in E_S}\left( X_i -\tilde{m}_i \right)   \right] \\
& \qquad = \mathbb{E}_{\mu_c ^{\sigma}} \left[\left(f(X)-\mathbb{E}_{\mu_c^{\sigma}} \left[f(X) \right] \right)\left(g(X)-\mathbb{E}_{\mu_c^{\sigma}} \left[g(X) \right] \right) \sum_{i \in E_S ^f}\left( X_i -\tilde{m}_i \right)   \right]  \label{e_second_taylor_1_1}\\
& \qquad \quad + \mathbb{E}_{\mu_c ^{\sigma}} \left[\left(f(X)-\mathbb{E}_{\mu_c^{\sigma}} \left[f(X) \right] \right)\left(g(X)-\mathbb{E}_{\mu_c^{\sigma}} \left[g(X) \right] \right) \sum_{i \in E_S ^g}\left( X_i -\tilde{m}_i \right)   \right] \label{e_second_taylor_1_2}
\end{align}

Then it follows that (see estimations of~\eqref{e_function_sec_2} for example)
\begin{align}
\left| T_{\eqref{e_second_taylor_1_1}} \right| &\leq \sum_{d \geq \frac{1}{2} d_{f,g} }\left| \mathbb{E}_{\mu_c ^{\sigma}} \left[\left(f(X)-\mathbb{E}_{\mu_c^{\sigma}} \left[f(X) \right] \right| \right)\left(g(X)-\mathbb{E}_{\mu_c^{\sigma}} \left[g(X) \right] \right) \sum_{\substack{i \in E_S^f, \\ \text{dist}(i, \supp f ) = d}}\left( X_i -\tilde{m}_i \right)   \right]  \\
& =\sum_{d \geq \frac{1}{2} d_{f,g} } \left| \cov_{\mu_c ^{\sigma}} \left( f(X), \left(g(X)-\mathbb{E}_{\mu_c^{\sigma}} \left[g(X) \right] \right) \sum_{\substack{i \in E_S^f, \\ \text{dist}(i, \supp f ) = d}}\left( X_i - \tilde{m}_i \right) \right) \right| \\
& \overset{Theorem~\ref{p_decay_of_correlations_gce}}{\lesssim} \sum_{d \geq \frac{1}{2} d_{f,g} } \| \nabla f \|_{L^2 (\mu_c ^{\sigma})} \| \nabla g \|_{L^4 (\mu_c ^{\sigma})} |\{ i \in E_S ^f \ | \ \text{dist}(i, \supp f ) = d \}|^{\frac{1}{2}} \exp \left( -Cd \right) \\
& \lesssim \| \nabla f \|_{L^2 (\mu_c ^{\sigma})} \| \nabla g \|_{L^4 (\mu_c ^{\sigma})}  |\supp f |^{\frac{1}{2}} \sum_{d \geq \frac{1}{2} d_{f,g} } \exp \left( -Cd \right) \\
& \lesssim \| \nabla f \|_{\infty} \| \nabla g \|_{\infty}  |\supp f |^{\frac{1}{2}} \exp \left( -C d_{f,g} \right).
\end{align}

Similarly, one gets
\begin{align}
\left| T_{\eqref{e_second_taylor_1_2}} \right| &\lesssim \| \nabla f \|_{\infty} \| \nabla g \|_{\infty}  |\supp g |^{\frac{1}{2}}\exp \left( -C d_{f,g} \right),
\end{align}
and thus 
\begin{align}
\left|T_{\eqref{e_second_taylor_1}} \right| &\lesssim \| \nabla f \|_{\infty} \| \nabla g \|_{\infty}  \left( |\supp f |^{\frac{1}{2}} +|\supp g |^{\frac{1}{2}} \right)    \exp \left( -C d_{f,g} \right) \frac{|\xi|}{\sqrt{N}} \\
& \leq \| \nabla f \|_{\infty} \| \nabla g \|_{\infty}\exp \left( -C d_{f,g} \right) |\xi|.
\end{align}

\medskip

Let us turn to the estimation of~\eqref{e_second_taylor_2}. It holds that
\begin{align}
& \mathbb{E}_{\mu_c ^{\sigma}} \left[\left(f(X)-\mathbb{E}_{\mu_c^{\sigma}} \left[f(X) \right] \right)\left(g(X)-\mathbb{E}_{\mu_c^{\sigma}} \left[g(X) \right] \right) \left(\sum_{i \in E_S} \left( X_i -\tilde{m}_i \right) \right)^2  \right]\\
&  = \cov_{\mu_c ^{\sigma}} \left( \left(f(X)-\mathbb{E}_{\mu_c^{\sigma}} \left[f(X) \right] \right)\left(g(X)-\mathbb{E}_{\mu_c^{\sigma}} \left[g(X) \right] \right)  ,   \left(\sum_{i \in E_S} \left( X_i -\tilde{m}_i \right) \right)^2    \right) \label{e_second_taylor_2_1} \\
&  \quad + \mathbb{E}_{\mu_c ^{\sigma}} \left[ \left(f(X)-\mathbb{E}_{\mu_c^{\sigma}} \left[f(X) \right] \right)\left(g(X)-\mathbb{E}_{\mu_c^{\sigma}} \left[g(X) \right] \right)\right] \mathbb{E}_{\mu_c ^{\sigma}} \left[\left(\sum_{i \in E_S} \left( X_i -\tilde{m}_i \right) \right)^2 \right]. \label{e_second_taylor_2_2}
\end{align}
Then Lemma~\ref{l_function_second_estimate} implies that
\begin{align}
\left| T_{\eqref{e_second_taylor_2_1}}\right| &\lesssim \| \nabla \left(  \left(f(X)-\mathbb{E}_{\mu_c^{\sigma}} \left[f(X) \right] \right)\left(g(X)-\mathbb{E}_{\mu_c^{\sigma}} \left[g(X) \right] \right) \right)  \|_{L^4 (\mu_c ^{\sigma})}|S| \\
& \leq \left( \| \nabla f \|_{\infty} \| g(X)-\mathbb{E}_{\mu_c^{\sigma}} \left[g(X) \right] \|_{L^4 (\mu_c ^{\sigma})} + \| \nabla g \|_{\infty} \| f(X)-\mathbb{E}_{\mu_c^{\sigma}} \left[f(X) \right] \|_{L^4 (\mu_c ^{\sigma})} \right) |S| \\
& \overset{Lemma~\ref{l_function_moment_estimate}}{\lesssim} \| \nabla f \|_{\infty}\| \nabla g \|_{\infty} |S|.
\end{align}
The second term~\eqref{e_second_taylor_2_2} is estimated via Theorem~\ref{p_decay_of_correlations_gce} and Lemma~\ref{l_variance_estimate} as follows:
\begin{align}
\left| T_{\eqref{e_second_taylor_2_2}}\right| &= \left|\cov_{\mu_c ^{\sigma}}\left( f(X)  ,  g(X) \right) \right| \mathbb{E}_{\mu_c ^{\sigma}} \left[\left(\sum_{i \in E_S} \left( X_i -\tilde{m}_i \right) \right)^2 \right] \\
& \lesssim \| \nabla f \|_{\infty}\| \nabla g \|_{\infty} \exp \left( -C d_{f,g} \right) |S|.
\end{align}
Thus plugging the estimates for~$T_{\eqref{e_second_taylor_2_1}}$ and~$T_{\eqref{e_second_taylor_2_2}}$ into~\eqref{e_second_taylor_2} yields
\begin{align}
\left|T_{\eqref{e_second_taylor_2}}\right| &\lesssim \| \nabla f \|_{\infty} \| \nabla g \|_{\infty} \frac{|S|}{N} \xi^2 \\
& \quad + \| \nabla f \|_{\infty}\| \nabla g \|_{\infty} \exp \left( -C d_{f,g} \right) \xi^2,
\end{align}
where we used~$|S| \leq N$ in the second term. \\

Lastly, we address the cubic term~\eqref{e_second_taylor_3}. H\"older's inequality followed by Lemma~\ref{l_fourth_moment} gives
\begin{align}
\left| T_{\eqref{e_second_taylor_3}} \right| &\lesssim \|\left(f(X)-\mathbb{E}_{\mu_c^{\sigma}} \left[f(X) \right] \right)\left(g(X)-\mathbb{E}_{\mu_c^{\sigma}} \left[g(X) \right] \right) \|_{L^4 (\mu_c ^{\sigma})} \| \sum_{i \in E_S} \left( X_i -\tilde{m}_i \right) \|_{L^4 (\mu_c ^{\sigma})}^3 \frac{|\xi|^3}{N^{\frac{3}{2}}} \\
& \lesssim \|\left(f(X)-\mathbb{E}_{\mu_c^{\sigma}} \left[f(X) \right] \right)\left(g(X)-\mathbb{E}_{\mu_c^{\sigma}} \left[g(X) \right] \right) \|_{L^4 (\mu_c ^{\sigma})} \frac{|E_S|^{\frac{3}{2}}}{N^{\frac{3}{2}}} |\xi|^3.
\end{align}
A combination of Schwarz inequality and Lemma~\ref{l_function_moment_estimate} implies
\begin{align}
\|\left(f(X)-\mathbb{E}_{\mu_c^{\sigma}} \left[f(X) \right] \right)\left(g(X)-\mathbb{E}_{\mu_c^{\sigma}} \left[g(X) \right] \right) \|_{L^4 (\mu_c ^{\sigma})} \leq \| \nabla f \|_{\infty} \| \nabla g \|_{\infty}.
\end{align}
Because~$|E_S| \leq 2\left(|\supp f |+ |\supp g | \right)  M \ln N$, it holds for~$N$ large enough that
\begin{align}
\frac{|E_S|^{\frac{3}{2}}}{N^{\frac{3}{2}}} \lesssim \frac{|\supp f | + |\supp g|}{N}.
\end{align}
Therefore we conclude
\begin{align}
\left| T_{\eqref{e_second_taylor_3}} \right| \lesssim \| \nabla f \|_{\infty} \| \nabla g \|_{\infty} \frac{|\supp f | + |\supp g|}{N} |\xi|^3.
\end{align}

Lastly, we sum up all the estimates we have obtained so far. That is,
\begin{align}
\left|T_{\eqref{e_second_taylor_original}} \right| &\leq  \left|T_{\eqref{e_second_taylor_0}} \right|+\left|T_{\eqref{e_second_taylor_1}} \right|+\left|T_{\eqref{e_second_taylor_2}} \right|+\left|T_{\eqref{e_second_taylor_3}} \right| \\
& \lesssim \| \nabla f \|_{\infty} \| \nabla g \|_{\infty} \left( \frac{|\supp f | + |\supp g | }{N}+  \exp\left(-Cd_{f,g}\right)  \right)\left( 1 + |\xi| + \xi^2 + |\xi|^3 \right).
\end{align}
Putting this estimate and Lemma~\ref{l_exchanging_exponential_terms} into~\eqref{e_second_der_further_decomposition_1} implies
\begin{align}
\left| T_{\eqref{e_cov_decomp_1}} \right| &\lesssim \| \nabla f \|_{\infty} \| \nabla g \|_{\infty} \left( \frac{|\supp f | + |\supp g | }{N}+  \exp\left(-Cd_{f,g}\right)  \right) \\
& \qquad \times \left( 1 + |\xi| + \xi^2 + |\xi|^3 \right) \exp \left(-C \xi^2 \right) \\
& \lesssim \| \nabla f \|_{\infty} \| \nabla g \|_{\infty} \left( \frac{|\supp f | + |\supp g | }{N}+  \exp\left(-Cd_{f,g}\right)  \right) \exp \left( -C \xi^2 \right).
\end{align}
\qed

\medskip

\begin{lemma} \label{l_2nd_der_inner_2}
It holds that
\begin{align}
\left|T_{\eqref{e_cov_decomp_2}} \right|, \left|T_{\eqref{e_cov_decomp_3}} \right|, \left|T_{\eqref{e_cov_decomp_4}} \right| \lesssim \| \nabla f \|_{\infty} \| \nabla g \|_{\infty} \frac{ |\supp f |+ |\supp g |}{N} \exp \left(-C\xi^2 \right).
\end{align}
\end{lemma}

\medskip 

\noindent \emph{Proof of Lemma~\ref{l_2nd_der_inner_2}.} \ Similar to the proof of Lemma~\ref{l_2nd_der_inner_1} we decompose~$T_{\eqref{e_cov_decomp_2}}$ as follows:
\begin{align}
T_{\eqref{e_cov_decomp_2}} & = \mathbb{E}_{\mu^{\sigma}} \left[ \exp \left( i \frac{1}{\sqrt{N}} \sum_{j \in F_S ^1} \left( X_j-m_j \right)\xi + i \frac{1}{\sqrt{N}} \sum_{i \in E_S} \left( \tilde{m}_i - m_i \right) \xi \right)  \right. \\
& \qquad \qquad  \times \mathbb{E}_{\mu_c ^{\sigma}} \left[ \exp \left( i \frac{1}{\sqrt{N}} \sum_{i \in E_S} \left( X_i -\tilde{m}_i \right) \xi \right) \left(g(X)-\mathbb{E}_{\mu_c^{\sigma}} \left[g(X) \right] \right)  \right] \\
& \qquad \qquad \qquad   \times  \left(\mathbb{E}_{\mu_c^{\sigma}} \left[f(X) \right]-\mathbb{E}_{\mu^{\sigma}} \left[f(X) \right] \right)   \\
& \qquad \qquad \qquad \quad \quad  \times \left. \mathbb{E}_{\mu^{\sigma}} \left[ \left. \exp \left( i \frac{1}{\sqrt{N}} \sum_{k \in F_S ^2} \left( X_j-m_j \right)\xi \right) \ \right| \  X_i, i \in E_S \cup F_S ^1  \right] \ \right],
\end{align}
First of all, Lemma~\ref{l_function_moment_estimate} implies
\begin{align}
\left| \mathbb{E}_{\mu_c ^{\sigma}} \left[ \exp \left( i \frac{1}{\sqrt{N}} \sum_{i \in E_S} \left( X_i -\tilde{m}_i \right) \xi \right) \left(g(X)-\mathbb{E}_{\mu_c^{\sigma}} \left[g(X) \right] \right)  \right] \right|&  \leq \mathbb{E}_{\mu_c ^{\sigma}} \left[ \left|g(X)-\mathbb{E}_{\mu_c^{\sigma}} \left[g(X) \right] \right|  \right] \\
& \lesssim \| \nabla g \|_{\infty}.
\end{align}
Next, we apply Lemma~\ref{l_conditional_function_difference} and Lemma~\ref{l_exchanging_exponential_terms} to obtain
\begin{align}
\left|T_{\eqref{e_cov_decomp_2}}\right| &\lesssim \| \nabla f \|_{\infty} \| \nabla g \|_{\infty}\exp \left( -C M \ln N\right)  \exp \left(-C\xi^2 \right)\\
& \qquad \times \mathbb{E}_{\mu^{\sigma}} \left[  \int \left( \sum_{\substack{i \in E_S , j \in F_S \\ |i-j| \leq R}} M_{ij}^2(y_j -z_j )^2 \right)^{\frac{1}{2}} \bar{\mu}^{\sigma} (dz)  \right] \\
&\overset{\eqref{e_mij_proof}}{\lesssim}   \| \nabla f \|_{\infty} \| \nabla g \|_{\infty}\exp \left( -C M \ln N\right)  \left( |\supp f| + |\supp g | \right)^{\frac{1}{2}} \exp \left(-C\xi^2 \right) \\
& \lesssim  \| \nabla f \|_{\infty} \| \nabla g \|_{\infty} \frac{ |\supp f |+ |\supp g |}{N} \exp \left(-C\xi^2 \right).
\end{align}
Similar calculations also imply
\begin{align}
\left|T_{\eqref{e_cov_decomp_3}}\right|, \left|T_{\eqref{e_cov_decomp_4}}\right| \lesssim \| \nabla f \|_{\infty} \| \nabla g \|_{\infty} \frac{ |\supp f |+ |\supp g |}{N} \exp \left(-C\xi^2 \right).
\end{align}
\qed

Now we are ready to give a proof of Lemma~\ref{l_2nd_der_inner}, which finishes the proof. \\

\noindent \emph{Proof of Lemma~\ref{l_2nd_der_inner}.} \ By lemma~\ref{l_2nd_der_inner_1} and Lemma~\ref{l_2nd_der_inner_2} it holds that
\begin{align}
&\left| \int_{\{  \left| \left(1/\sqrt{N}\right) \xi \right| \leq \delta \}} \mathbb{E}_{\mu^{\sigma}} \left[ \exp \left( i \frac{1}{\sqrt{N}} \sum_{k=1}^{N} \left( X_k - m_k \right) \xi \right) \right. \right. \\
& \qquad \qquad \qquad \qquad \qquad \qquad \left. \left. \times \vphantom{\exp \left( i \frac{1}{\sqrt{N}} \sum_{k=1}^{N} \left( X_k - m_k \right) \xi \right)} \left( f(X) - \mathbb{E}_{\mu^{\sigma}} \left[f(X) \right] \right)\left( g(X) - \mathbb{E}_{\mu^{\sigma}} \left[g(X) \right] \right) \right] d\xi \right| \\
& \qquad \lesssim \int_{\mathbb{R}} \| \nabla f \|_{\infty} \| \nabla g \|_{\infty} \left( \frac{|\supp f | + |\supp g | }{N}+  \exp\left(-Cd_{f,g}\right)  \right) \exp \left( -C \xi^2 \right) d\xi \\
& \qquad \lesssim \| \nabla f \|_{\infty} \| \nabla g \|_{\infty} \left( \frac{|\supp f | + |\supp g | }{N}+  \exp\left(-Cd_{f,g}\right)  \right).
\end{align}
\qed 

\section{Appendix : Proof of Lemma~\ref{l_inverse_fourier}}

\medskip 
Let~$\zeta, \eta : \mathbb{R}^N \to \mathbb{R}$ be given. Let us fix~$\sigma = \sigma(m)$ and~$m = m(\sigma)$ as in Section~\ref{s_setting_and_main_results}. We introduce auxiliary external fields with~$\zeta, \eta$ in the definition of gce and ce. More precisely, let us denote for~$\alpha, \beta \in \mathbb{R}$,
\begin{align}
\mu^{\sigma, \alpha, \beta} (dx) &:= \frac{1}{Z} \exp \left( \sigma \sum_{i=1}^{N} x_i + \alpha \zeta(x) + \beta \eta (x) - H(x) \right) dx, \\
\mu_m ^{\alpha, \beta} (dx) &:= \frac{1}{Z} \mathds{1}_{\left\{ \frac{1}{N} \sum_{i=1}^{N} x_i = m \right\} }  \exp \left(  \sigma \sum_{i=1}^{N} x_i + \alpha \zeta(x) + \beta \eta(x) - H(x) \right) \mathcal{L}^{N-1} (dx). 
\end{align}
In particular, one should observe that
\begin{align} \label{e_relation_auxiliary_external_fields}
\left. \mu^{\sigma, \alpha, \beta} (dx) \right|_{\alpha, \beta = 0} = \mu^{\sigma} (dx), \qquad
\left. \mu_m ^{ \alpha, \beta} (dx) \right|_{\alpha, \beta = 0} = \mu_m (dx).
\end{align}
Let us consider associated free energies~$A_{gce}^{\zeta, \eta}$,~$A_{ce}^{\zeta, \eta}$ defined by
\begin{align}
A_{gce}^{\zeta, \eta} (\alpha, \beta) &:= \frac{1}{N} \ln \int \exp \left( \sigma \sum_{i=1}^{N} x_i + \alpha \zeta (x) + \beta \eta (x) - H(x) \right) dx , \\
A_{ce}^{\zeta, \eta} (\alpha, \beta) &:= \frac{1}{N} \ln \int_{ \left\{ \frac{1}{N} \sum_{i=1}^N x_i = m  \right\} } \exp \left( \sigma \sum_{i=1}^{N} x_i + \alpha \zeta (x) + \beta \eta (x) - H(x) \right) \mathcal{L}^{N-1} (dx).
\end{align}

First, a straightforward calculation yields

\begin{align}
\left.  \frac{d}{d \alpha} A_{gce}^{\zeta, \eta} \right|_{\alpha, \beta = 0} &= \frac{1}{N} \mathbb{E}_{\mu^{\sigma}} \left[ \zeta (X) \right], \qquad \left.  \frac{d}{d \beta} A_{gce}^{\zeta, \eta} \right|_{\alpha, \beta = 0} = \frac{1}{N} \mathbb{E}_{\mu^{\sigma}} \left[ \eta (X) \right], \\
\left.  \frac{d}{d \alpha} A_{ce}^{\zeta, \eta} \right|_{\alpha, \beta = 0} &= \frac{1}{N} \mathbb{E}_{\mu_m} \left[ \zeta (X) \right], \qquad \left.  \frac{d}{d \beta} A_{ce}^{\zeta, \eta} \right|_{\alpha, \beta = 0} = \frac{1}{N} \mathbb{E}_{\mu_m} \left[ \eta (X) \right], \\
\left.  \frac{d^2}{d \alpha d \beta} A_{gce}^{\zeta, \eta} \right|_{\alpha, \beta = 0} &= \frac{1}{N} \cov_{\mu^{\sigma}} \left( \zeta (X) , \eta (X) \right), \quad 
\left.  \frac{d^2}{d \alpha d \beta} A_{ce}^{\zeta, \eta} \right|_{\alpha, \beta = 0} = \frac{1}{N} \cov_{\mu_m} \left( \zeta (X) , \eta (X) \right).  \label{e_derivatives_auxiliary_gce_ce}
\end{align}

Next, the Cramer's representation yields
\begin{align}
&A_{ce}^{\zeta, \eta} (\alpha, \beta) - A_{gce}^{\zeta,  \eta} (\alpha, \beta) \\
& \qquad = \frac{1}{N} \ln \frac{\int_{ \left\{ \frac{1}{N} \sum_{i=1}^{N} x_i = m \right\} } \exp \left( \sigma \sum_{i=1}^{N} x_i + \alpha \zeta (x) + \beta \eta(x) - H(x)  \right) \mathcal{L}^{N-1} (dx)  }{\int_{\mathbb{R}} \exp \left( \sigma \sum_{i=1}^{N} x_i + \alpha \zeta (x) + \beta \eta (x) - H(x) \right) dx  } \\
& \qquad = \frac{1}{N} \ln \frac{\int_{ \left\{ \frac{1}{\sqrt{N}} \sum_{i=1}^{N} (x_i - m) = 0 \right\} } \exp \left( \sigma \sum_{i=1}^{N} x_i + \alpha \zeta (x) + \beta \eta(x) - H(x)  \right) \mathcal{L}^{N-1} (dx)  }{\int_{\mathbb{R}} \exp \left( \sigma \sum_{i=1}^{N} x_i + \alpha \zeta (x) + \beta \eta (x) - H(x) \right) dx  }. \label{e_aux_cramer_rep}
\end{align}

Let~$W = (W_1, \cdots, W_N)$ be a real-valued random vector distributed according to~$\mu^{\sigma, \alpha, \beta}$ and~$g_{\alpha, \beta}$ be the density of random variable
\begin{align}
\frac{1}{\sqrt{N}} \sum_{i=1}^{N} (W_i - m).
\end{align}
Then in~\eqref{e_aux_cramer_rep}, it holds that
\begin{align} \label{e_aux_density}
A_{ce}^{\zeta, \eta} (\alpha, \beta) - A_{gce}^{\zeta \eta} (\alpha, \beta) = \frac{1}{N} \ln g_{\alpha, \beta}(0).
\end{align}
Note also that an application of inverse Fourier transformation yields
\begin{align}
2 \pi g_{\alpha, \beta} (0) = \int_{\mathbb{R}} \mathbb{E}_{\mu^{\sigma, \alpha, \beta}} \left[ \exp \left( i \frac{1}{\sqrt{N}} \sum_{i=1}^{N} \left( W_i - m \right) \xi \right) \right] d \xi  .
\end{align}

Again, a direct calculation with a help of identity~\eqref{e_relation_auxiliary_external_fields} implies
\begin{align}
\left. 2 \pi g_{\alpha, \beta} (0) \right|_{\alpha, \beta = 0} &= \int_{\mathbb{R}} \mathbb{E}_{\mu^{\sigma}} \left[\exp \left( i \frac{1}{\sqrt{N}} \sum_{i=1}^{N} \left( X_i - m \right) \xi \right) \right] dx, \\
\left. 2 \pi \frac{d}{d\alpha} g_{\alpha, \beta} (0) \right|_{\alpha, \beta = 0} &= \int_{\mathbb{R}} \mathbb{E}_{\mu^{\sigma}} \left[ \left( \zeta (X) - \mathbb{E}_{\mu^{\sigma}} \left[ \zeta (X) \right] \right)  \exp \left( i \frac{1}{\sqrt{N}} \sum_{i=1}^{N} \left( X_i - m \right) \xi \right) \right] dx, \\
\left. 2 \pi \frac{d}{d\beta} g_{\alpha, \beta} (0) \right|_{\alpha, \beta = 0} &= \int_{\mathbb{R}} \mathbb{E}_{\mu^{\sigma}} \left[ \left( \eta (X) - \mathbb{E}_{\mu^{\sigma}} \left[ \eta (X) \right] \right)  \exp \left( i \frac{1}{\sqrt{N}} \sum_{i=1}^{N} \left( X_i - m \right) \xi \right) \right] dx, \\
\left. 2 \pi \frac{d^2}{d\alpha d\beta } g_{\alpha, \beta} (0) \right|_{\alpha, \beta = 0} &= \int_{\mathbb{R}} \mathbb{E}_{\mu^{\sigma}} \left[  \left( \zeta (X) - \mathbb{E}_{\mu^{\sigma}} \left[ \zeta (X) \right] \right) \left( \eta (X) - \mathbb{E}_{\mu^{\sigma}} \left[ \eta (X) \right] \right) \vphantom{\exp \left( i \frac{1}{\sqrt{N}} \sum_{i=1}^{N} \left( X_i - m \right) \xi \right)}  \right. \\
& \qquad \qquad \qquad \qquad \left. \times \exp \left( i \frac{1}{\sqrt{N}} \sum_{i=1}^{N} \left( X_i - m \right) \xi \right) \right] dx. \label{e_derivatives_density_auxiliary}
\end{align}
Then a combination of~\eqref{e_derivatives_auxiliary_gce_ce},~\eqref{e_aux_density} and~\eqref{e_derivatives_density_auxiliary} yields the desired result.
\qed

\section*{Acknowledgment}
This research has been partially supported by NSF grant DMS-1712632. The authors are thankful to many people discussing the problem and helping to improve the preprint. Among them are Tim Austin, Frank Barthe, Marek Biskup, Pietro Caputo, Jean-Dominique Deuschel, Max Fathi, Andrew Krieger, Michel Ledoux, Sangchul Lee, Thomas Liggett, Guido Mont\'ufar, Felix Otto, Andr\'e Schlichting, Daniel Ueltschi, and Tianqi Wu. The authors want to thank Marek Biskup, UCLA and KFAS for financial support.

\bibliographystyle{alpha}
\bibliography{bib}

\end{document}